%% file: paper.tex
\documentclass{article}

\input{ryantibs}

\def\d{\textup{d}}  
\def\one{\mathds{1}}
\def\iidsim{\overset{\mathrm{iid}}{\sim}} 
 
\DeclareMathOperator{\infconv}{\texttt{\textup{\#}}}
\DeclareMathOperator{\epi}{epi} 
\DeclareMathOperator{\prox}{prox}
\DeclareMathOperator{\diam}{diam}    
\DeclareMathOperator{\interior}{int}
\DeclareMathOperator{\softmax}{softmax} 

\newcommand{\ie}{i.e.,\ }
\newcommand{\eg}{e.g.,\ }  

\usepackage{xcolor}
\usepackage{setspace}
\usepackage{pdfpages}
\usepackage{etoolbox}
\usepackage[shortcuts]{extdash}

\usepackage{xcolor}
\definecolor{TypalBlue}{HTML}{7086DF}  
\definecolor{TypalRed}{HTML}{D32B08}  
\definecolor{TypalGreen}{HTML}{3cb371}  
\newcommand{\rottext}[1]{\begin{minipage}[t]{0.2in}\rotatebox{90}{\parbox[c]{1.1in}{\hspace*{10pt}\centering #1}}\end{minipage}\hspace*{-8pt}}

\usepackage{tikz}

\title{Laplace Meets Moreau: Smooth Approximation to Infimal Convolutions Using
  Laplace's Method}     
\author{Ryan J.\ Tibshirani$^a$ \and Samy Wu Fung$^b$ \and Howard Heaton$^c$
  \and Stanley Osher$^d$}   
\date{$^a$University of California, Berkeley \quad
$^b$Colorado School of Mines \quad
$^c$Typal Academy \quad
$^d$University of California, Los Angeles
}

\begin{document}
\maketitle

\begin{abstract}
We study approximations to the Moreau envelope---and infimal convolutions more
broadly---based on Laplace's method, a classical tool in analysis which ties
certain integrals to suprema of their integrands. We believe the connection
between Laplace's method and infimal convolutions is generally deserving of more
attention in the study of optimization and partial differential equations, since
it bears numerous potentially important applications, from proximal-type
algorithms to solving Halmiton-Jacobi equations.   
\end{abstract}

\section{Introduction} 

Infimal convolutions are of core importance in mathematical optimization and
partial differential equations (PDEs). The most well-known special case of an
infimal convolution is the Moreau envelope, due to \citet{moreau1962fonctions,
  moreau1965proximite}, which (along with its counterpart, the proximal
operator) is a key tool in convex and variational analysis, and in numerical
algorithms for optimization. More broadly, beyond the Moreau envelope, infimal
convolutions appear as solutions in a class of Hamilton-Jacobi equations in
PDEs. 



Laplace's method, due to \citet{laplace1774memoire}, is a tool for approximating 
integrals that finds applications in many areas of mathematics, statistics,
physics, and computer science. More specifically, it provides a way to precisely
approximate an integral whose integrand becomes increasingly peaked around its
maximum value. To researchers in statistics and machine learning, Laplace's
method is perhaps most familiar from its use in Bayesian inference, where it
leads to an approximation of the posterior distribution in terms of a Gaussian
distribution centered at the maximum a posteriori (MAP) estimate. 

These two ideas are actually closely connected: Laplace's method provides a
natural way to approximate an infimal convolution. This has been noted and used
(albeit somewhat indirectly) by some authors in the past; see Section
\ref{sec:related_work}. However, we believe the connection between infimal
convolutions and Laplace's method is not as widely appreciated as it should be,
especially as it relates to sampling, which is a way to view (and numerically
approximate) the integrals that appear in Laplace's method. Thus, the current
paper places the connection between infimal convolutions and Laplace's method
front and center. It is not really our intent to claim novelty in developing or
formalizing this connection. Instead, our goal is to highlight both its elegance 
and utility in the hope that it gains better recognition, and potentially, sees
further applications. 

\paragraph{Outline.} 

In what follows, we first review preliminary concepts and related work. In
Section \ref{sec:smooth_approximation}, we describe the use of Laplace's method
to approximate infimal convolutions, and we give interpretations from various
perspectives. In Section \ref{sec:asymptotic_theory}, we derive approximation
guarantees. In Section \ref{sec:monte_carlo} we cover sampling techniques, and
in Section \ref{sec:applications_examples}, we walk through applications in
optimization and PDEs, with illustrative examples. 



\subsection{Infimal convolution}

\begin{figure}[t]
\centering 
\includegraphics{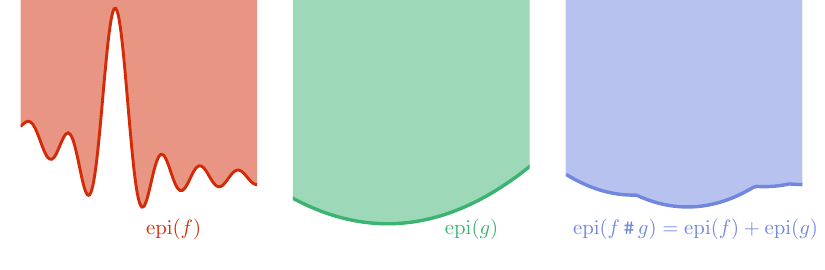} 
\caption{\small 
  Illustration of infimal convolution as epigraph addition
  \eqref{eq:epi_addition} (here $g = \|\cdot\|_2^2 / 2$).} 
\label{fig:epi_addition}
\end{figure}

Let $f, g : \R^d \to \R$ be arbitrary real-valued functions. The \emph{infimal
  convolution} (or simply the inf convolution) of $f$ and $g$ is another
function, denoted $f \infconv g : \R^d \to \R$, which is defined by       
\begin{equation} 
\label{eq:inf_conv}
(f \infconv g)(x) = \inf_y \, \Big\{ f(y) + g(x-y) \Big\}.
\end{equation}
The notion of an infimal convolution originated with \citet{fenchel1951convex},
while a related idea was independently developed over a series of papers by
\citet{bellman1961new, bellman1962mathematical, bellman1962maximum,
  bellman1963maximum}. Influenced by Fenchel, the papers by
\citet{moreau1963inf} and \citet{rockafellar1963convex} serve as the basis for 
what is now considered the modern definition and treatment of inf convolutions,
with Moreau's work providing the name ``inf convolution'', as well. For more on
infimal convolutions, we refer to \citet{moreau1970inf,stromberg1994infimal, 
  rockafellar2009variational}. The latter book introduces the notation
$\infconv$ for the inf convolution operator that we use in this paper, which is
meant to remind the reader of the addition operator, because (as Rockafellar and
Wets emphasize) infimal convolution acts as addition on the space of epigraphs:      
\begin{equation}
\label{eq:epi_addition}
\epi(f \infconv g) = \epi(f) + \epi(g),
\end{equation}
as long as the infimum defining $(f \infconv g)(x)$ is attained whenever
finite. Here $\epi(f) = \{ (x,t) : f(x) \leq t\}$ is the epigraph of $f$,
similarly for $\epi(g)$, and $A+B = \{a+b : a \in A, \, b \in B\}$ is the
usual (Minkowski) sum of sets $A,B$. Figure \ref{fig:epi_addition} gives an
illustration. 

A particularly important special case of an inf convolution is the Moreau
envelope. This plays a central role in various aspects of optimization, from
theoretical to practical, and is covered in the next subsection. Inf
convolutions also play an important role as solutions to certain Hamilton-Jacobi
equations. To elaborate, let $H : \R^d \to \R$ be a convex Hamiltonian, 
and consider the first-order PDE  
\begin{equation}
\begin{alignedat}{2}
\label{eq:hj_prob}
\partial_t u + H(\nabla u) &= 0, && \quad t>0, \\
u(x,0) &= f(x), && \quad t=0,
\end{alignedat}
\end{equation}
where $\partial_t u$ denotes the derivative of $u$ with respect to $t$, and 
$\nabla u$ denotes its gradient with respect to $x$. By the Hopf-Lax formula,
(\eg Theorem 6 in Chapter 3.3.2 of \citet{evans2010partial}), the solution is     
\begin{equation}
\label{eq:hj_sol}
u(x,t) = \inf_y \, \Big\{ f(y) + t H^*\Big( \frac{x-y}{t} \Big) \Big\}, 
\quad t>0,
\end{equation}
where $H^*$ is the conjugate (also called the Legendre-Fenchel transform) of
$H$. In other words, the solution to the Hamilton-Jacobi PDE \eqref{eq:hj_prob}
at time $t>0$ is given by the infimal convolution $u(x,t) = (f \infconv \,
tH^*(\cdot/t))(x)$. 

We finish this subsection with a useful general fact about infimal
convolutions. Fix any point $x$, assume that there is a unique point $y_x$ which  
attains the infimum in \eqref{eq:inf_conv}, and assume $g$ is differentiable on 
$\R^d$. Then under some additional regularity conditions 
(\eg Theorem 10.13 and Corollary 10.14 of \citet{rockafellar2009variational}),
the inf convolution $f \infconv g$ is differentiable at $x$ and       
\begin{equation}
\label{eq:inf_conv_grad}
\nabla (f \infconv g)(x) = \nabla g(x-y_x).
\end{equation}
We will return to this formula shortly, in the case of the Moreau envelope,
discussed next.      

\subsection{Moreau regularization} 

\begin{figure}[t]
\centering 
\includegraphics{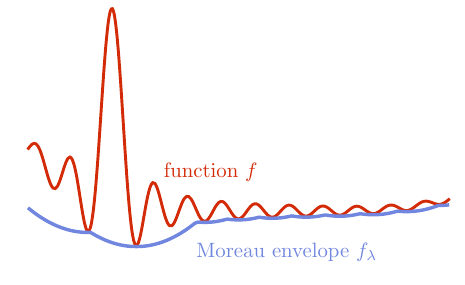} 
\caption{\small
  Illustration of the Moreau envelope $f_\lambda = f \infconv \|\cdot\|_2^2/
  (2\lambda)$ (for a particular $\lambda$).} 
\label{fig:moreau_env}
\end{figure}

For $g = \|\cdot\|_2^2/ (2\lambda)$, with $\lambda>0$ a fixed constant, the
infimal convolution $f \infconv g$ is called the \emph{Moreau envelope} of $f$
at the level $\lambda$, named after the pioneering work of
\citet{moreau1962fonctions, moreau1965proximite}, and is denoted by 
\begin{equation}
\label{eq:moreau_env}
f_\lambda(x) = \inf_y \, \Big\{ f(y) + \frac{1}{2\lambda} \|x-y\|_2^2 \Big\}.  
\end{equation}
Figure \ref{fig:moreau_env} gives an illustration. Intimately connected to this
is the \emph{proximal operator} of $\lambda f$, denoted by     
\begin{equation}
\label{eq:prox}
\prox_{\lambda f}(x) = \argmin_y \, \Big\{ f(y) + \frac{1}{2\lambda} \|x-y\|_2^2
\Big\}, 
\end{equation}
where we use \smash{$\argmin_y F(y)$} to denote the set of minimizers of a
function $F$, reducing to \smash{$\argmin_y F(y) = z$} if the set of minimizers
is a singleton $\{z\}$. Importantly, if $f$ is convex, then the operator
\smash{$\prox_{\lambda f}$} is guaranteed to be single-valued (rather than
set-valued): it maps each input $x$ to a unique point \smash{$\prox_{\lambda
    f}(x)$}.        
      
The Moreau envelope \eqref{eq:moreau_env} and its associated proximal operator
\eqref{eq:prox} are ubiquitous throughout convex and variational analysis, as
well as optimization. In terms of theory, proximal operators admit various
important connections to subdifferentials, conjugates, and monotone operators;
see, \eg \citet{rockafellar2009variational, bauschke2011convex}. In terms of
algorithms, proximal operators serve as a building block for a number of
operator splitting techniques for nonsmooth, constrained, large-scale
optimization, which includes forward-backward splitting, Douglas-Rachford
splitting, and the alternating direction method of multipliers (ADMM); see, \eg
\citet{boyd2011distributed, combettes2011proximal, parikh2013proximal,
  beck2017first, ryu2022large}.     

An early and influential contribution on the algorithmic front is called the
\emph{proximal point algorithm}, due to \citet{rockafellar1976monotone}. Given a 
function $f$ to be minimized, we fix $\lambda>0$, initalize $x_0 \in \R^d$, and
then repeat the iterations, for $k=1,2,3,\dots$:
\begin{equation}
\label{eq:prox_point}
x_k = \prox_{\lambda f}(x_{k-1}).
\end{equation}
Returning to the gradient formula for $f \infconv g$ in
\eqref{eq:inf_conv_grad}, for a closed convex function $f$ (and with $g = 
\|\cdot\|_2^2/(2\lambda)$, which together satisfy the conditions needed for this
gradient formula), we have 
\begin{equation}
\label{eq:moreau_grad}
(\nabla f_\lambda)(x) = \frac{x - \prox_{\lambda f}(x)}{\lambda} 
\iff \prox_{\lambda f}(x) = x - \lambda (\nabla f_\lambda)(x).
\end{equation}
This is true regardless of the smoothness of $f$. Applying
\eqref{eq:moreau_grad} to \eqref{eq:prox_point}, we see that Rockafellar's
proximal point iteration can be rewritten as  
\begin{equation}
\label{eq:moreau_grad_desc}
x_k = x_{k-1} - \lambda (\nabla f_\lambda)(x_{k-1}),
\end{equation}
which is the same as gradient descent on the Moreau envelope. 

\subsection{Laplace's method} 

Let $\varphi : \R^d \to \R$ be twice continuously differentiable and $h : \R^d
\to \R$ be continuous. \emph{Laplace's method} (also called Laplace
approximation) provides an asymptotic equivalence for an integral that becomes 
increasingly peaked around the global minimizer $x^\star$ of $\varphi$, assumed
to be unique:
\begin{equation}
\label{eq:laplace_asymp}
\int h(x) \exp(-t\varphi(x)) \, \d{x} \;\sim\; \frac{(2\pi/t)^{d/2}} 
{\exp(t\varphi(x^\star))} \frac{h(x^\star)} {\sqrt{\det(\nabla^2
    \varphi(x^\star))}},
\end{equation}
where we use $a(t) \sim b(t)$ as shorthand for $a(t)/b(t) \to 1$ as $t \to
\infty$.  

In short, applications of Laplace's method can be found throughout statistics,
probability, and machine learning (not to mention its uses in mathematics and
physics). Examples include 
Bayesian computation and inference (\eg \citet{kass1991laplace}),   
higher-order asymptotics (\eg \citet{shun1995laplace}),
mean-field theory (where Laplace's method is often called the \emph{saddle 
  point approximation}, \eg \citet{mezard2009information}, Chapter 2), 
Gaussian processes (\eg \citet{rasmussen2006gaussian}, Chapter 3),
and Bayesian deep learning (\eg \citet{daxberger2021laplace}).

Now we present a formal statement of the validity of Laplace's approximation    
\eqref{eq:laplace_asymp}. We allow $\varphi$ to have an arbitrary domain
$\cK$, and reparametrize by $t = 1/\delta$ (taking $\delta \to 0^+$ in the
asymptotic limit), since it will be more convenient for our purposes later.      

\newcommand{\thmLaplaceMethod}[1]{
Let $\varphi : \cK \to \R$ be continuous over a compact set $\cK \subseteq  
\R^d$. Assume that $\varphi$ has a unique global minimizer $x^\star$ in the
interior of $\cK$, and $\varphi$ is twice continuously differentiable on a
neighborhood of $x^\star$, with strictly positive definite Hessian $\nabla^2
\varphi(x^\star)$. Then for any continuous function $h : \cK \to \R$,  
\begin{equation}
\label{eq:laplace_method_#1}
\sqrt{\det(\nabla^2 \varphi(x^\star))} \;\; \frac{\exp(\varphi(x^\star)/\delta)}    
{(2\pi\delta)^{d/2}} \int_\cK h(x) \exp(-\varphi(x) / \delta) \, \d{x} \;\to\;
h(x^\star), \quad \text{as $\delta \to 0^+$}.     
\end{equation}
This conclusion extends to the case where $\cK$ is not compact (\eg $\cK = 
\R^d$) provided there exists $\epsilon>0$ such that the sublevel set
$\cS_\epsilon = \{ x \in \cK : \varphi(x) \leq \varphi(x^\star) + \epsilon\}$
is bounded and \smash{$\int_\cK |h(x)| \exp(-d \varphi(x) / (2\epsilon)) \,
  \d{x} < \infty$}.  
}

\begin{theorem}
\label{thm:laplace_method} 
\thmLaplaceMethod{main}
\end{theorem}

Proofs of the asymptotic convergence of Laplace's method can be found in any
standard reference on the topic. For completeness, we provide a proof of Theorem
\ref{thm:laplace_method} in Appendix \ref{app:laplace_method}, based on the
simple and elegant arguments given in \citet{bach2021approximating}.    

The form of the approximation given in Theorem \ref{thm:laplace_method} (or
equivalently in \eqref{eq:laplace_asymp}) is well-suited for traditional
applications of Laplace's method. In such problems, we seek to avoid computing
an integral in $\exp(-t\varphi)$, the left-hand side in
\eqref{eq:laplace_asymp}, and approximate it using a minimum of $\varphi$ (or
maximum of $-\varphi$), the right-hand side in \eqref{eq:laplace_asymp}. For
example, in Bayesian inference, this method can be used to approximate the
posterior distribution using a Gaussian centered at the maximum a posteriori
(MAP) estimate. Here, the sample size $n$ plays the role of $t$ in
\eqref{eq:laplace_asymp}, making the approximation more accurate for larger
sample sizes.   

For our purposes however, a different form of the Laplace approximation will be
more convenient. This is because our motivation is really the \emph{opposite} of
the traditional one: we seek to avoid computing a minimum of $\varphi$, and we 
instead approximate it using an integral involving $\exp(-t\varphi)$. In our
setting, we would not want the approximation to feature normalizing constants
such as $\exp(t\varphi(x^\star))$ and \smash{$\sqrt{\det(\nabla^2
    \varphi(x^\star))}$}, since they are unknown (recall, $\varphi(x^\star)$ is
the minimum that we are trying to approximate in the first place). Fortunately,
we can use the following \emph{self-normalized} version of Laplace's method:    
\begin{equation}
\label{eq:laplace_asymp_sn}
\frac{\int_\cK h(x) \exp(-t\varphi(x)) \, \d{x}}{\int_\cK \exp(-t\varphi(x)) \,
  \d{x}} \;\sim\; h(x^\star).
\end{equation}
This follows directly from \eqref{eq:laplace_method_main}, by using the latter
to approximate separately the integrals in the numerator and denominator in
\eqref{eq:laplace_asymp_sn} (the common factor \smash{$(2\pi/t)^{d/2} /
  (\exp(t\varphi(x^\star)) \sqrt{\det(\nabla^2\varphi(x^\star))} )$} cancels
out). Because this just relies on two applications of Laplace's method, Theorem
\ref{thm:laplace_method} gives the formal conditions under which
\eqref{eq:laplace_asymp_sn} is valid. Later in Section
\ref{sec:asymptotic_theory}, we will give a generalization of this result: 
instead of requiring $\varphi$ to be twice differentiable on a neighborhood of
its minimizer, we only require it to be locally H{\"o}lder continous.  

\subsection{Related work}
\label{sec:related_work}

The literature related to the topic of our paper is vast, and below is our
attempt to give a broad overview of the various related areas, without giving a  
comprehensive review of any one in particular.

\paragraph{Viscosity solutions in PDEs.}  

In a certain sense, the use of Laplace's method in order to approximate
solutions in Hamilton-Jacobi equations dates back to seminal work on  
\emph{viscosity solutions} by \citet{crandall1983viscosity,
  crandall1984properties} (see also \citet{evans1980solving}). A canonical way
to construct viscosity solutions in Hamilton-Jacobi equations is the called the
\emph{vanishing viscosity} method, where we solve a modified PDE that has an
additional diffusion term with a small viscosity parameter, and then send the
viscosity parameter to zero. Interestingly, as we review in Section 
\ref{sec:viscous_burgers}, for the Hamilton-Jacobi equation \eqref{eq:hj_prob} 
where $H = \|\cdot\|_2^2 / 2$ (also called Burgers' equation), the use of a
Laplacian diffusion term leads to exactly the same approximation as Laplace's
method applied to the original solution. As far as we understand, some but not
all authors in the PDE literature on viscosity make the connection to Laplace's
approximation explicit. For example, it is not discussed in the early papers in  
the 1980s by Crandall, Lions, and Evans, but it can be found in 
\citet{evans2010partial} in Chapter 4.5.2. We have not yet seen the use of
Laplace's approximation for general Hamilton-Jacobi equations \eqref{eq:hj_prob}
(for general $H$, as we propose in this paper), and as we discuss at the end of
Section \ref{sec:viscous_burgers}, we are unsure as to whether there is a
viscosity-like representation for such an approximation in general. 

In the total variation-based image denoising literature, the posterior mean
formula \eqref{eq:laplace_prox_exp}, a special case of the Laplace approximation
\eqref{eq:laplace_inf_prox_exp} for infimal convolutions that we consider in
this paper, has been studied by \citet{louchet2008modeles,
  louchet2013posterior}. This was extended by \citet{darbon2021bayesian,
  darbon2021connecting}, who derive rigorous approximations guarantees using  
connections to Hamilton-Jacobi PDEs. Another line of work that draws connections 
between viscous Hamilton-Jacobi PDEs, proximal operators, and Moreau envelopes
was initiated by \citet{chaudhari2018deep} and further developed by
\citet{heaton2023global, osher2023hamilton}. As we explain in Sections 
\ref{sec:exponential_tilting} and \ref{sec:viscous_burgers}, in this paper we
arrive at the identical approximation for the proximal operator as that given in
\citet{darbon2021connecting, heaton2023global, osher2023hamilton}, albeit from a
different perspective: by directly applying Laplace's method, instead of relying
on viscosity. By casting the approximation through the lens of Laplace's method,
we are able to seamlessly extend it to handle arbitrary infimal convolutions. We
are also able to make less stringent assumptions (on the functions in question)
for the approximation theory that we derive in Section
\ref{sec:asymptotic_theory}.    

\paragraph{Sampling and optimization.}

As we explore in Sections \ref{sec:prox_point} and \ref{sec:bregman_prox},
Laplace's approximation \eqref{eq:laplace_prox_exp} of the proximal operator
(and \eqref{eq:laplace_inf_prox_exp} for inf convolutions more generally) leads
to various sampling-based methods for optimization. These methods are 
\emph{zeroth-order}: they depend only on function evaluations (and not
gradients, as would be the case in first-order methods). The connections between
sampling and optimization are quite deep (in both directions---using sampling to    
optimize, and using optimization to sample). Classical examples include
stochastic gradient descent \citep{robbins1951stochastic}, simulated annealing
\citep{kirkpatrick1983optimization}, and Langevin dynamics
\citep{welling2011bayesian}. More recently, stochastic localization and
diffusion models have taken a center stage in machine learning (\eg
\citet{sohldickstein2015deep, song2019generative, ho2020denoising,
  song2021score, elalaoui2022sampling, montanari2023sampling}),   
and remain an extremely active topic of research.  

In the setting of zeroth-order optimization in particular, the connections to 
sampling are also rich, dating back to \citet{matya1965random}. Related to  
the Laplace approximation of the proximal point algorithm (studied in Section 
\ref{sec:prox_point}) is the idea of Gaussian smoothing from
\citet{nesterov2017random}, who study a gradient-free optimization algorithm 
which approximates a directional derivative with a finite difference scheme
based on a random pertubation of the parameter. The motivation and focus in
their work, as with much of the literature in zeroth-order optimization, is
quite different than ours. We return to this discussion in Section
\ref{sec:prox_point}.      

\section{Smooth approximation}
\label{sec:smooth_approximation}

Let $f,g : \R^d \to \R$ be continuous functions. Define for any fixed $x \in 
\R^d$ and $\delta>0$:
\begin{equation} 
\label{eq:laplace_inf_prox}
y_x^\delta = \frac{\int y \exp\Big( \frac{-f(y) - g(x-y)}{\delta} \Big) \, 
  \d{y}}{\int \exp\Big( \frac{-f(y) - g(x-y)}{\delta} \Big) \, \d{y}}. 
\end{equation}
Notice that each coordinate \smash{$(y_x^\delta)_i$} in
\eqref{eq:laplace_inf_prox} is the self-normalized Laplace approximation from 
\eqref{eq:laplace_asymp_sn}, applied to $\varphi_x(y) = f(y) + g(x-y)$, with
$h(y) = y_i$. In other words, we can view \smash{$y_x^\delta$} as approximating
the minimizer of the inf convolution criterion at $x$, 
\[
  y_x^\delta \approx \argmin_y \, \big\{ f(y) + g(x-y) \big\}.
\]
This approximation becomes exact as $\delta \to 0^+$ under mild conditions on 
$f,g$, as the theory in the next section will make precise (see Corollary
\ref{cor:laplace_inf_conv}). It is important to note that convexity of $f,g$ is 
not required. Below, we discuss interpretations of the approximation
\eqref{eq:laplace_inf_prox}, from different perspectives.      

\subsection{Exponential tilting}
\label{sec:exponential_tilting}

Define a density by
\begin{equation}
\label{eq:p_delta_g}
p^\delta_{g,x}(y) = \frac{\exp(-g(x-y) / \delta)}{\int \exp(-g(x-y) / \delta) \, 
\d{y}}.
\end{equation}
We can rewrite \eqref{eq:laplace_inf_prox} as 
\begin{equation} 
\label{eq:laplace_inf_prox_exp}
y_x^\delta = \frac{\E_{Y \sim p^\delta_{g,x}} [ Y \exp(-f(Y) / \delta) ]} 
{\E_{Y \sim p^\delta_{g,x}} [ \exp(-f(Y) / \delta) ]}.
\end{equation}
This can be interpreted as the expectation in a model in which we first sample
$Y$ from the density \smash{$p^\delta_{g,x}$}, which (assuming that $g$ is
minimized at the origin) is centered at $x$ and increasingly peaked for smaller
$\delta$, and then \emph{exponentially tilt} by $-f/\delta$, which upweights
the samples that lead to smaller values of $f$. Note that we can also interpret 
\eqref{eq:laplace_inf_prox_exp} from the Bayesian perspective: consider a 
Bayesian model with likelihood \smash{$p^\delta_{g,x} \propto e^{-g(x-\cdot)}$},
and prior $\propto e^{-f}$. In this context, the quantity \smash{$y_x^\delta$}
represents the \emph{posterior mean}.    

In the special case of $g = \|\cdot\|_2^2 / (2\lambda)$, observe that
\smash{$p^\delta_{g,x}$} in \eqref{eq:p_delta_g} is the $N(x, \delta \lambda I)$ 
density, and \eqref{eq:laplace_inf_prox_exp} leads to
\begin{equation} 
\label{eq:laplace_prox_exp}
y_x^\delta = \frac{\E_{Y \sim N(x, \delta \lambda I)}[ Y \exp(-f(Y) / \delta) ]}
{\E_{Y \sim N(x, \delta \lambda I)} [ \exp(-f(Y) / \delta) ]},    
\end{equation}
which recovers the proximal approximation formula in \citet{darbon2021bayesian,
  osher2023hamilton}, who arrived at this result from a different perspective,
as explained next.  

\subsection{Viscous Burgers' equation}
\label{sec:viscous_burgers}

As discussed previously, the Hopf-Lax formula \eqref{eq:hj_sol} gives the 
solution to the Hamilton-Jacobi PDE \eqref{eq:hj_prob}. When $H =
\|\cdot\|_2^2 / 2$, problem \eqref{eq:hj_prob} reduces to what is
known as Burgers' equation, and Laplace's approximation can be understood from
the perspective of what is called \emph{viscosity} in the PDE literature. In
particular, consider the viscous Burgers' equation (Example 2 in Chapter 4.5.2
of \citet{evans2010partial}),       
\begin{equation}
\begin{alignedat}{2}
\label{eq:viscous_burgers_prob}
\partial_t u^\delta + \frac{1}{2} \|\nabla u^\delta\|_2^2 &=
\frac{\delta}{2} \Delta u^\delta, && \quad t>0, \\  
u^\delta(x,0) &= f(x), && \quad t=0,
\end{alignedat}
\end{equation}
where $\Delta u$ is the Laplacian of $u$ with respect to $x$. By using a change 
of variables \smash{$v^\delta(x,t) = \exp(-u^\delta(x,t) / \delta)$} (also known 
as the Cole-Hopf transform), problem \eqref{eq:viscous_burgers_prob} becomes the  
heat equation, with the initial condition \smash{$v^\delta(x,0) = \exp(-f(x) / 
  \delta)$}. We can use the fundamental solution of the heat equation, and
translate back to our original parametrization, to yield the solution
\begin{equation}
\label{eq:viscous_burgers_sol}
u^\delta(x,t) = - \delta \log\Bigg( \frac{1}{(2\pi\delta t)^{d/2}} \int
  \exp\bigg(\frac{-f(y) - \|x-y\|_2^2 / (2t)}{\delta} \bigg) \, \d{y} \Bigg).
\end{equation}
This approximates the solution $u(x,t)$ in \eqref{eq:hj_sol}, \ie it
approximates the Moreau envelope $f_t$ of $f$ (since, recall, $H = \|\cdot\|_2^2 
/ 2$). A seminal result by \citet{crandall1984two} is that \smash{$u^\delta(x,t)
  \to u(x,t)$} as $\delta \to 0^+$, uniformly over all $x \in \R^d$ and all
compact intervals of time $t \geq 0$. Further results are available in
\citet{darbon2021bayesian, darbon2021connecting, heaton2023global,
  osher2023hamilton}. 

As observed by the aforementioned authors, we can combine
\eqref{eq:viscous_burgers_sol} with the Moreau gradient formula
\eqref{eq:moreau_grad} to obtain an approximation to the proximal  
map. First we differentiate \eqref{eq:viscous_burgers_sol} with respect to $x$,
\[
\nabla u^\delta(x,t) = \frac{\int \frac{x-y}{t} \exp\Big(\frac{-f(y) - 
    \|x-y\|_2^2 / (2t)}{\delta} \Big) \, \d{y}} {\int \exp\Big(\frac{-f(y) -
    \|x-y\|_2^2 / (2t)}{\delta} \Big) \, \d{y}}.
\]  
Then based on \eqref{eq:moreau_grad}, we approximate \smash{$\prox_{tf}(x)$}
using $x - t \nabla u^\delta(x,t)$. Observe that  
\[
x - t \nabla u^\delta(x,t) = \frac{\int y \exp\Big(\frac{-f(y) - \|x-y\|_2^2 /
    (2t)}{\delta} \Big) \, \d{y}} {\int \exp\Big(\frac{-f(y) - \|x-y\|_2^2 /
    (2t)}{\delta} \Big) \, \d{y}}
\]
is precisely the self-normalized Laplace approximation in 
\eqref{eq:laplace_inf_prox} with $g = \|\cdot\|_2^2 / (2t)$. Recall, this also
has the equivalent form \eqref{eq:laplace_prox_exp}, expressed in terms of
expectations with respect to $N(x, \delta t I)$. 

The fact that Laplace's approximation to the proximal operator can be 
alternatively derived via viscosity in Burgers' equation is quite interesting
(and to emphasize once again, this was the path taken by previous work to arrive 
at the same formula for the proximal approximation that we get from Laplace's 
method). This begs the question: for general $g$, is there such a viscosity-like 
representation for \eqref{eq:laplace_inf_prox}? That is, can we find a
viscosity-like modification to the Hamilton-Jacobi PDE \eqref{eq:hj_prob}, with
$H^* = g$, whose solution leads to the formula \eqref{eq:laplace_inf_prox}? 
While viscosity is studied for general Hamilton-Jacobi equations (\eg Chapter
10 of \citet{evans2010partial}), as far as we can tell, the viscous
Hamilton-Jacobi equation using a Laplacian diffusion term does not lead to an
solution that coincides with Laplace's approximation of the inf convolution in
general, in the way that it does when $H = \|\cdot\|_2^2 / 2$. Investigating
whether we can express \eqref{eq:laplace_inf_prox} in terms of a viscosity-like
pertubation of the Hamilton-Jacobi equation \eqref{eq:hj_prob} may be an
interesting direction for future investigation.    
  
   


\subsection{Smoothed set projection}

\begin{figure}[t]    
\centering   
\hspace*{-10pt}
\includegraphics[width=0.35\textwidth, page=2]{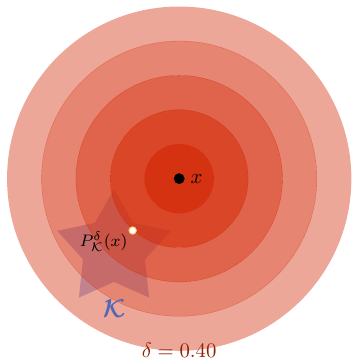} 
\hspace*{-15pt}
\includegraphics[width=0.35\textwidth, page=4]{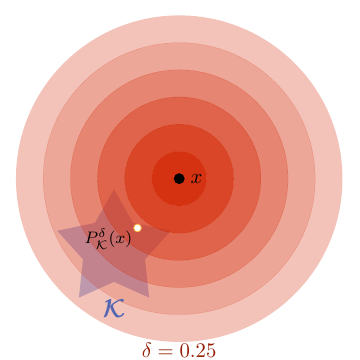} 
\hspace*{-15pt}
\includegraphics[width=0.35\textwidth, page=6]{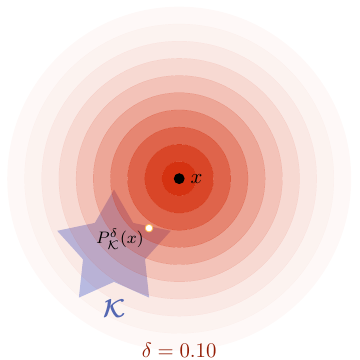}  
\hspace*{-10pt} 
\caption{\small
  Examples of the Laplace approximation \smash{$P_{\cK}^\delta(x)$} to the
  projection $P_\cK(x)$ of a point $x$ onto a set $\cK$. In each panel 
  \smash{$P_{\cK}^\delta(x)$} is denoted by a white dot, and is defined by   
  \eqref{eq:laplace_proj} for a particular value of $\delta$. As $\delta$ 
  decreases (from left to right), the conditional density of $Y \,|\, Y   \in
  \cK$ becomes increasingly peaked around $P_\cK(x)$. Each outer red circle
  about $x$ represents a successive standard deviation \smash{$\sqrt\delta$}
  in the contours of the sampling distribution $N(x, \delta I)$.}   
\label{fig:laplace_proj}
\end{figure}

Returning to the posterior mean formula \eqref{eq:laplace_prox_exp} when 
$g = \|\cdot\|_2^2 / 2$, consider taking $f = I_\cK$, the characteristic
function of a set $\cK \subseteq \R^d$, 
\[
I_\cK(x) = 
\begin{cases}
0 & \text{if $x \in \cK$} \\   
\infty & \text{otherwise}.
\end{cases}  
\]
In this case, the quantity being approximated is the proximal operator
\smash{$y_x = \prox_{I_\cK}(x)$} of $I_\cK$ evaluated at $x$, \ie the projection 
$y_x = P_\cK(x)$ of $x$ onto the set $\cK$, concretely
\[
P_\cK(x) = \argmin_{y \in \cK} \, \|x-y\|_2^2.
\]
Introducing the notation \smash{$P_\cK^\delta(x) = y_x^\delta$} for the 
approximation in \eqref{eq:laplace_prox_exp}, note that this simplifies to   
\begin{equation}
\label{eq:laplace_proj}
P_\cK^\delta(x) = \E_{Y \sim N(x, \delta I)}[ Y \,|\, Y \in \cK], 
\end{equation}
This formula is highly intuitive: we take an average according to a certain
density over $\cK$, which for small $\delta$, ends up being nearly flat at points
far away from $P_\cK(x)$, and more peaked close to $P_\cK(x)$. Figure
\ref{fig:laplace_proj} gives an illustration.  
 
If $\cK$ is convex then the approximation \eqref{eq:laplace_proj} has the
property that \smash{$P_\cK^\delta(x) \in \cK$} for any $\delta>0$; for
nonconvex $\cK$, this no longer needs to be true, and we could have
\smash{$P_\cK^\delta(x)$} lying outside of $\cK$. However, convexity is not
required in order to guarantee \smash{$P_\cK^\delta(x) \to P_\cK(x)$} as $\delta
\to 0^+$. We only require a mild condition on the boundary of $\cK$ in a
neighborhood of $P_\cK(x)$ (see Corollary \ref{cor:laplace_proj}).

\subsection{Integral convolution}

Lastly, we make the simple observation that in the current setting the Laplace
approximation brings an inf convolution to an ordinary (integral)
convolution. Generalizing from \eqref{eq:laplace_inf_prox}, suppose we are
interested in   
\[
h(y_x) = h\Big( \argmin_y \, \big\{ f(y) + g(x-y) \big\} \Big),
\]
for a given function $h : \R^d \to \R$ (we also assume that the argmin above is
unique). Then the corresponding self-normalized Laplace approximation from
\eqref{eq:laplace_asymp_sn} is  
\begin{equation}
\label{eq:laplace_h_inf_prox}
\frac{\int h(y) \exp\Big( \frac{-f(y)-g(x-y)}{\delta} \Big) \, \d{y}}{\int 
\exp\Big( \frac{-f(y)-g(x-y)}{\delta} \Big) \, \d{y}}. 
\end{equation}
Assuming $g$ is an even function (\ie $g(x) = g(-x)$), this can be expressed as   
\begin{equation}
\label{eq:laplace_h_int_conv}
\bigg( \frac{(h e^{-f/\delta}) * e^{-g/\delta}} {e^{-f/\delta} * e^{-g/\delta}}
\bigg) (x),
\end{equation}
where $(u * v)(x) = \int u(y) v(y-x) \, \d{y} = \int u(y) v(x-y) \, \d{y}$, the
last equality holding for even $v$. As convolution is generally understood as a
smoothing operation, this formulation informally lends support to the idea that 
Laplace method's provides a smooth approximation to $h(y_x)$.     

More formally, we can use properties of convolutions to infer about the
smoothness of \eqref{eq:laplace_h_inf_prox} as a function of $x$. If $u,v$ are 
integrable and $v$ has integrable partial derivatives, then a standard
fact (which can be verified using Fourier transforms) is that 
\[
\frac{\partial(u * v)}{\partial x_i} = u * \frac{\partial v}{\partial x_i},
\quad i = 1,\dots,n.
\]
Applying this to \eqref{eq:laplace_h_int_conv} (and using the quotient rule and
chain rule as needed), we see that the approximation is differentiable as many
times as $g$ is (\ie it is infinitely differentiable for a choice such as $g =
\|\cdot\|_2^2 / 2$).     

\section{Asymptotic theory}
\label{sec:asymptotic_theory}

We analyze the asymptotic validity ($\delta \to 0^+$) of the self-normalized
version of Laplace's approximation, in a way that generalizes what is known
classically \eqref{eq:laplace_asymp_sn}, which requires twice differentiability
(recall Theorem \ref{thm:laplace_method}). We break our presentation in what
follows into two parts, depending on whether the minimizer in question lies in
the interior of its domain.   
  
\subsection{Minimizer in the interior}

First, we study the Laplace approximation of a function $\varphi$ whose
minimizer is in the interior of its domain. The proof of the next theorem is
given in Appendix \ref{app:laplace_sn}.

\newcommand{\thmLaplaceSelfNormalized}[1]{
Let $\varphi : \cK \to \R$ be a continuous function, and assume that $\varphi$ 
has a unique global minimizer $x^\star$ that lies in the interior of $\cK$, 
admits a bounded sublevel set $\cS_\epsilon = \{ x \in \cK : \varphi(x) \leq
\varphi(x^\star) + \epsilon\}$ for some $\epsilon>0$, and satisfies
\smash{$\int_\cK \exp(-\varphi(x) / \epsilon) \, \d{x} < \infty$}. Assume that
on a neighborhood $U$ of $x^\star$, we have  
\begin{equation}
\label{eq:local_holder_#1}
\varphi(x) - \varphi(x^\star) \leq a \|x-x^\star\|_2^q, \quad \text{for all $x 
  \in \cK \cap U$},   
\end{equation}
for constants $a,q > 0$. Then for any continuous $h : \cK \to \R$ with
\smash{$\int_\cK |h(x)| \exp(-d \varphi(x) / (q\epsilon)) \, \d{x} < \infty$},        
\begin{equation}
\label{eq:laplace_sn_#1}
\frac{\int_\cK h(x) \exp(-\varphi(x) / \delta) \, \d{x}}{\int_\cK
  \exp(-\varphi(x) / \delta) \, \d{x}} \;\to\; h(x^\star), \quad \text{as
  $\delta \to 0^+$}.  
\end{equation}
}

\begin{theorem}
\label{thm:laplace_sn} 
\thmLaplaceSelfNormalized{main}
\end{theorem}

The condition \eqref{eq:local_holder_main} is essentially a type of local
H{\"o}lder continuity assumption on $\varphi$, on a neighborhood $U$ of its 
minimizer $x^\star$, with an arbitrary exponent $q>0$.\footnote{We say 
  ``essentially'' because, technically, H{\"o}lder continuity is stronger, and
  would require $\varphi(x) - \varphi(x^\star) \leq a \|x-x^\star\|_2^q$ for 
  all pairs $x,y \in U$, whereas the condition \eqref{eq:local_holder_main}
  only requires this inequality to hold at pairs $x,x^\star \in U$.} 
This is considerably weaker than assuming that $\varphi$ is twice differentiable
on a neighborhood of $x^\star$, as in Theorem \ref{thm:laplace_method}. (We note
that if $\varphi$ is locally twice continuously differentiable, then it
satisfies \eqref{eq:local_holder_main} with $q=2$, which can be verified using a 
Taylor expansion. Meanwhile, if $\varphi$ is only locally continuously
differentiable, then it satisfies \eqref{eq:local_holder_main} with $q=1$, as 
$\|\nabla \varphi\|_2$ has a finite maximum on any compact subset containing its
minimizer.) 

\begin{figure}[t]  
\centering  
\small
\begin{tabular}{ccc}
$\varphi(x) = \frac{9}{40} + \frac{(x-1)^4}{20} + \frac{\sin(10\pi x)}{40x}$ & 
$\varphi(x) = |x|^{0.5}$ & 
$\varphi(x) = \displaystyle \sum_{k=0}^{100} 0.3^k \cos(23^k \pi x) $ \\[10pt]
\includegraphics[page=1]{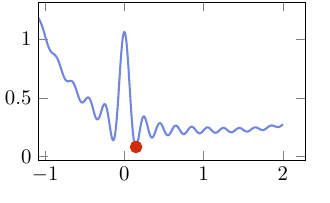} 
\hspace*{-10pt} &
\hspace*{-10pt}
\includegraphics[page=2]{images/1D_functions/1D-function-plots.pdf} 
\hspace*{-10pt} &
\hspace*{-10pt}
\includegraphics[page=3]{images/1D_functions/1D-function-plots.pdf}
\hspace*{-10pt} 
\end{tabular}
\caption{\small
  Examples of the Laplace approximation to the minimizer $x^\star$ of
  $\varphi$. The panels potray different functions $\varphi$, each of which
  satisfies the local H{\"o}lder assumption \eqref{eq:local_holder_main}, but is
  nonconvex and nonsmooth in the traditional sense. In each panel the Laplace
  approximation is denoted by a red dot, and computed using the left-hand side 
  of \eqref{eq:laplace_sn_main}, with $h$ being the identity map and
  $\delta=10^{-6}$.}   
\label{fig:laplace_sn}
\end{figure}

In order to be able to weaken the assumption from local twice differentiability
to local H{\"o}lder continuity, the self-normalized aspect of the
approximation in \eqref{eq:laplace_sn_main} is key, because in general
the explicit normalizing factors in a statement like
\eqref{eq:laplace_method_main} would require precise knowledge of the local
growth rate of $\varphi(x) - \varphi(x^\star)$. To be clear, in
\eqref{eq:local_holder_main}, we only need to know that this local growth rate
is a power of $\|x-x^\star\|_2$, without needing to know the exponent, in order
to compute the approximation \eqref{eq:laplace_sn_main}. Figure
\ref{fig:laplace_sn} gives a few illustrations.  

Next, we apply Theorem \ref{thm:laplace_sn} to an infimal convolution.         

\begin{corollary}
\label{cor:laplace_inf_conv}
Let $f,g : \R^d \to \R$ be continuous functions, and define $\varphi_x(y) = f(y)
+ g(x-y)$. Assume that $\varphi_x$ has a unique global minimizer $y_x$, has a
bounded sublevel set $\cS_\epsilon = \{ y \in \R^d : \varphi_x(y) \leq
\varphi_x(y_x) + \epsilon\}$ for some $\epsilon>0$, and satisfies \smash{$\int
  \exp(-\varphi_x(y) / \epsilon) \, \d{y} < \infty$}. Assume further that on a
neighborhood $U_x$ of $y_x$,     
\begin{equation}
\label{eq:local_holder}
\varphi_x(y) - \varphi_x(y_x) \leq a \|y-y_x\|_2^q, \quad \text{for all $y \in
  U_x$},   
\end{equation}
for constants $a,q > 0$. Then for any continuous $h : \R^d \to \R$ with
\smash{$\int |h(y)| \exp(-d\varphi_x(y) / (q\epsilon)) \, \d{y} < \infty$}, 
\[
\frac{\int h(y) \exp\Big( \frac{-f(y)-g(x-y)}{\delta} \Big) \, \d{y}}{\int 
\exp\Big( \frac{-f(y)-g(x-y)}{\delta} \Big) \, \d{y}} \;\to\; h(y_x), \quad
\text{as $\delta \to 0^+$}. 
\]
An important special case is given by taking $h(y) = y_i$, $i = 1,\dots,d$, for
which the above conclusion translates as follows: provided \smash{$\int
  \|y\|_\infty \exp((-f(y)-g(x-y)) / \epsilon) \, \d{y} < \infty$}, it holds for
\smash{$y_x^\delta$} as defined in \eqref{eq:laplace_inf_prox} that           
\[
y_x^\delta \to y_x, \quad \text{as $\delta \to 0^+$}. 
\vspace{4pt} 
\]
\end{corollary}

Once again we note that the local H{\"o}lder condition \eqref{eq:local_holder}
on $\varphi_x$ is quite weak. For example, this condition holds if $f,g$ are
each convex on a neighborhood of $y_x$ (since convex functions are locally
Lipschitz), or even if $f,g$ are weakly convex on a neighborhood of $y_x$ (see 
Appendix \ref{app:local_holder_weakly_convex} for a precise statement and
verification of this claim). While (weak) convexity is of central interest in
many applications of optimization, we emphasize that no kind of convexity is
required for \eqref{eq:local_holder}. Local H{\"o}lder continuity, with a
possibly fractional exponent $q$, can be satisfied by fairly exotic and highly
nonconvex functions (recall the functions from Figure
\ref{fig:laplace_sn}). 

Lastly, we note that with a more sophisticated analysis it is likely that the
local H{\"o}lder condition in \eqref{eq:local_holder} could be weakened further
into one on the local modulus of continuity, where we require that the modulus
of continuity have subexponential growth around zero (this is sufficient to
imply that the error term analyzed in step 3 of the proof vanishes; see
Appendix \ref{app:laplace_sn}). However, we do not pursue such an extension.  

\subsection{Minimizer on the boundary}

We consider the case in which the minimizer of $\varphi$ lies on the boundary of 
its domain. The proof of the next theorem is given in Appendix
\ref{app:laplace_sn_bdry}.        

\newcommand{\thmLaplaceSelfNormalizedBoundary}[1]{
Under the conditions of Theorem \ref{thm:laplace_sn}, assume instead that the
unique global minimizer $x^\star$ of $\varphi$ lies on the boundary of the
closed set $\cK$. Furthermore, assume that $\cK$ is full-dimensional and
star-shaped, in a local sense around $x^\star$: precisely, writing $B(x,r) = \{
y \in \R^d : \|y-x\|_2 \leq r \}$ for the ball of radius $r$ centered at $x$, we 
assume that there exists $r_0>0$ such that $\cK \cap B(x^\star,r_0)$ has
positive Lebesgue measure, and     
\begin{equation}
\label{eq:local_star_#1}
x \in \cK \cap B(x^\star,r_0) \implies \alpha x + (1-\alpha) x^\star \in \cK,
\quad \text{for all $\alpha \in [0,1]$}.
\end{equation}
Then the same conclusion as in \eqref{eq:laplace_sn_#1} holds.              
}

\begin{theorem}
\label{thm:laplace_sn_bdry} 
\thmLaplaceSelfNormalizedBoundary{main}
\end{theorem}

Next, we apply Theorem \ref{thm:laplace_sn_bdry} to a projection.      

\begin{corollary}
\label{cor:laplace_proj}
Let $\cK \subseteq \R^d$ and $x \in \R^d$ be such that $x$ has a unique
projection $y_x = P_\cK(x)$ onto the closed set $\cK$. Assume that there exists 
$r_0>0$ such that $\cK \cap B(y_x,r_0)$ has positive Lebesgue measure and
condition \eqref{eq:local_star_main} holds, with $y_x$ in place of
$x^\star$. Then for any continuous $h : \R^d \to \R$ with \smash{$\int_\cK
  |h(y)| \exp(-d\|x-y\|_2^2/2) \, \d{y} < \infty$},      
\[
\frac{\int_\cK h(y) \exp\Big( \frac{-\|x-y\|_2^2}{\delta} \Big) \, \d{y}}
{\int_\cK \exp\Big( \frac{-\|x-y\|_2^2}{\delta} \Big) \, \d{y}} \;\to\; h(y_x),
\quad \text{as $\delta \to 0^+$}. 
\]
An important special case is given by taking $h(y) = y_i$, $i = 1,\dots,d$, for
which the above conclusion translates as follows: provided \smash{$\int_\cK
  \|y\|_\infty \exp(-d\|x-y\|_2^2/2) \, \d{y} < \infty$}, it holds for
\smash{$P_\cK^\delta(x)$} as defined in \eqref{eq:laplace_proj} that          
\[
P_\cK^\delta(x) \to P_\cK(x), \quad \text{as $\delta \to 0^+$}. 
\]
\end{corollary}

The assumptions on the set $\cK$ used above, in the theorem and corollary, are
not strong. They are met if $\cK$ is a convex body (see Appendix
\ref{app:local_star_convex} for a precise statement and proof of this
claim). Yet the end results will continue to hold well outside of convexity.
Fundamentally, the stated assumptions on $\cK$ are really used as sufficient 
conditions to ensure that
\begin{equation}
\label{eq:local_hausdorff}
\cH^{d-1}(\cK \cap \partial B(x^\star,r)) \geq c r^{d-1}, \quad \text{for all $r 
  \leq r_0$}, 
\end{equation}
where $c>0$ is constant, and $\cH^{d-1}$ denotes Hausdorff measure of dimension
$d-1$. Informally, this condition says that $\cK$ should act in the way we would
expect of a ``regular'' full-dimensional set with nonempty interior, locally
around $x^\star$. Indeed, if $\cK = \R^d$, then condition
\eqref{eq:local_hausdorff} is met because $\cH^{d-1}(\partial B(x^\star,r)) = c
r^{d-1}$ for any $r$ and a constant $c>0$ depending only on $d$. Generally, one
can interpret \eqref{eq:local_hausdorff} as requiring the set $\cK$ to retain a
constant fraction of the surface measure over the boundary of the ball
$B(x^\star,r)$, for enough small $r$, which can still be satisfied by a highly
nonconvex set $\cK$, as long as $\cK$ has a somewhat ``regular'' behavior  
around its boundary (recall Figure \ref{fig:laplace_proj}). We note
that the condition \eqref{eq:local_hausdorff} could likely be weakened further
with a more sophisticated analysis, but we do not pursue this.

\section{Monte Carlo sampling}
\label{sec:monte_carlo}

In this section, we briefly discuss how to use sampling to approximate Laplace's
approximation \eqref{eq:laplace_prox_exp} of the proximal operator, and more
generally \eqref{eq:laplace_inf_prox_exp} for the minimizer of $f +
g(x-\cdot)$. We should note at the outset that Monte Carlo sampling is a rich
field with many powerful tools, and the ideas we describe here are only very 
basic. Likely, more advanced tools from the Monte Carlo literature
(and even from quasi-Monte Carlo) could be applied to improve accuracy and
efficiency.  

For the proximal formula \eqref{eq:laplace_prox_exp}, we can approximate this
using a sample average over Gaussian draws:        
\begin{gather}
\label{eq:gaussian_samples}
Y_i \iidsim N(x, \delta \lambda I), \; i = 1,\dots, N, \\
\label{eq:laplace_prox_avg}
y_x^{\delta,N} = \frac{\sum_{i=1}^N Y_i \exp(-f(Y_i) / \delta)}
{\sum_{i=1}^N \exp(-f(Y_i) / \delta)},    
\end{gather}
We note that \eqref{eq:laplace_prox_avg} can also be succinctly written as   
\begin{equation}
\label{eq:laplace_softmax}
y_x^{\delta,N} = (Y_1,\dots,Y_N)^\T \softmax(-f(Y_1)/\delta, \dots,
-f(Y_N)/\delta),  
\end{equation}
where for a vector $v = (v_1,\dots,v_N) \in \R^N$, we denote \smash{$\softmax(v) 
  = (e^{v_1} / \sum_{i=1}^N e^{v_i}, \dots, e^{v_N} / \sum_{i=1}^N
  e^{v_i})$}. Practical implementations of the softmax operator (such as that
in SciPy) commonly shift the exponents in order to avoid overflow, instead
computing \smash{$\softmax(v) = (e^{v_1-a} / \sum_{i=1}^N e^{v_i-a}, \dots, 
  e^{v_N-a} / \sum_{i=1}^N e^{v_i-a})$}, for a scalar value $a$. Any value
will contribute a common factor $e^a$ which cancels in the numerator and
denominator; but a careful choice such as \smash{$a = \max_{i=1,\dots,N} v_i$}
can lead to be better numerical accuracy (\eg see
\citet{blanchard2021accurately}). For this reason, it can be advantageous to
implement \eqref{eq:laplace_prox_avg} using \eqref{eq:laplace_softmax} in
practice (to take advantage of built-in shifting for numerical robustness) and
the same comment applies to all of the formulae involving exponentially-weighted
averages in the remainder of this section.     

For general $g$, we can approximate \eqref{eq:laplace_inf_prox_exp} analogously
using a sample average over suitable draws:  
\begin{gather}
\label{eq:p_delta_g_samples}
Y_i \iidsim p^\delta_{g,x}, \; i = 1,\dots, N, \\
\label{eq:laplace_inf_prox_avg}
y_x^{\delta,N} = \frac{\sum_{i=1}^N Y_i \exp(-f(Y_i) / \delta)}
{\sum_{i=1}^N \exp(-f(Y_i) / \delta)},    
\end{gather}
where \smash{$p^\delta_{g,x}$} is the density in \eqref{eq:p_delta_g}. Sampling from
\smash{$p^\delta_{g,x}$} will really only be possible in certain special
cases. For example, aside from the case $g = \|\cdot\|_2^2$, where
\smash{$p^\delta_{g,x}$} is a Gaussian density, we note that when $g =  
\|\cdot\|_1$, the density \smash{$p^\delta_{g,x}$} is a product of Laplace
densities (one for each coordinate).  Outside of such special cases, we can
use importance sampling, where instead of \eqref{eq:p_delta_g_samples},
\eqref{eq:laplace_inf_prox_avg}, we compute: 
\begin{gather}
\label{eq:q_samples}
Y_i \iidsim q, \; i = 1,\dots, N, \\
\label{eq:laplace_inf_prox_avg_iw}
y_x^{\delta,N} = \frac{\sum_{i=1}^N Y_i w(Y_i) \exp(-f(Y_i) / \delta)} 
{\sum_{i=1}^N w(Y_i) \exp(-f(Y_i) / \delta)},     
\end{gather}
where \smash{$w(y) = p^\delta_{g,x}(y) / q(y)$}, and $q$ is a user-chosen 
proposal density (whose support contains that of \smash{$p^\delta_{g,x}$}). In
general, a goal in choosing the proposal density $q$ would be to minimize 
the variance of the weighted sample average in
\eqref{eq:laplace_inf_prox_avg_iw}. This is a nontrivial task, but various
practical solutions have been developed by \citet{hesterberg1988advances,
  veach1997robust, owen2000safe}, among others.    

\section{Applications and examples}
\label{sec:applications_examples}

In what follows, we walk through applications of the Laplace approximations
proposed and studied above to problems in PDEs and optimization. To be clear, in
each case, we do not intend to produce or compete with state-of-the-art
solutions for the problem at hand. We only aim to demonstrate the broad 
applicability and portability of Laplace's approximation, through relatively
simple experiments. We focus on low-dimensional problems were sampling is
fairly easy (naive Monte Carlo or importance sampling works fairly
well). Higher-dimensional problems would call for more advanced sampling
techniques.    

\subsection{Hamilton-Jacobi equations}
\label{sec:hj_experiments}

We examine the use of Laplace's method for solving the Hamilton-Jacobi (HJ)  
equation \eqref{eq:hj_prob}. In particular, we study how $\delta$ and the number
of samples $N$ used in the approximation \eqref{eq:q_samples},
\eqref{eq:laplace_inf_prox_avg_iw} influence the quality of the (approximate)
HJ solution, for different convex Hamiltonians $H(x) = \|x\|_p^p / p$, $p \in 
\{1.1,2,5,10\}$, in different dimensions $d \in \{2,5,10\}$. We also set $f(x) =
\|x\|_1$.  

Recall that the exact HJ solution $u(x,t)$ is given in \eqref{eq:hj_sol}. 
Consider approximating this by first computing \smash{$y_x^\delta$} as in
\eqref{eq:laplace_inf_prox}, and then using \smash{$u^\delta(x,t) =
  f(y_x^\delta) + g(y_x^\delta - x)$}, with   
\[
g(x) = t H^*\bigg (\frac{x}{t} \bigg), \quad \text{and} \quad H^*(x) =
\frac{1}{q} \|x\|_q^q  
\]
for $1/q + 1/p = 1$. 
For any fixed $t$, we approximate \smash{$y_x^\delta$} with
\smash{$y_x^{\delta,N}$} in \eqref{eq:q_samples},
\eqref{eq:laplace_inf_prox_avg_iw}, where $N$ denotes the number of
samples. Then, we calculate     
\begin{equation}
\label{eq:lapprox_hj_sol}
u^{\delta,N} (x,t) = f(y_x^{\delta,N}) + g(y_x^{\delta,N} - x),
\end{equation}
at 1000 uniformly sampled values of $x \in [-10, 10]^d$ and $t \in [10^{-1},
1]$. The error of the approximation \smash{$u^{\delta,N}$} is measured by
calculating the magnitude of the HJ residual,   
\begin{equation}
\label{eq:lapprox_hj_res}
r(x,t) = \big| \partial_t u^{\delta,N}(x,t) + H(\nabla u^{\delta,N}(x,t)) \big|,
\end{equation}
and utlimately averaging this over the sampled values of $x$ and $t$.

To compute the Laplace approximation \smash{$y_x^{\delta,N}$} in
\eqref{eq:q_samples}, \eqref{eq:laplace_inf_prox_avg_iw} we choose the proposal 
density $q$ to be uniform over $[-10,10]^d$. Practically, limiting the domain in
this way has little effect (integrating over larger domains lead to similar
results). Given \smash{$u^{\delta, N}$} in \eqref{eq:lapprox_hj_sol}, we then
compute the HJ residual \eqref{eq:lapprox_hj_res} using PyTorch's automatic
differentiation functionality \citep{paszke2019pytorch} for the associated
derivatives.          

\begin{figure}[t]
\centering 
\hspace*{-5pt}
\rottext{\qquad\small Residual} 
\hspace*{-10pt}
\subfloat[Dimension $d=2$, norm $p=2$]{
\includegraphics[width=0.33\textwidth, page=4]{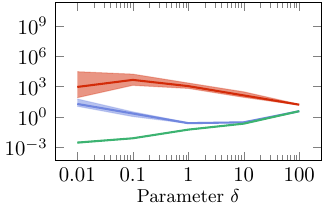}} 
\hspace*{-10pt}
\subfloat[Dimension $d=2$, norm $p=5$]{
\includegraphics[width=0.33\textwidth, page=7]{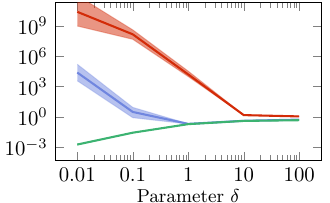}}
\hspace*{-10pt}
\subfloat[Dimension $d=2$, norm $p=10$]{
\includegraphics[width=0.33\textwidth, page=10]{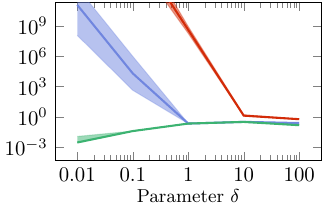}} 
\caption{\small
  Residuals for HJ approximation in dimension $d=2$, for $H = \|\cdot\|_p^p / 
  p$, with $p \in \{2,5,10\}$. Each panel shows \textcolor{TypalRed}{$N = 10$
    samples via red}; \textcolor{TypalBlue}{$N=10^3$ samples via blue}; and 
  \textcolor{TypalGreen}{$N=10^5$ samples via green}. We can see that as the
  dimension grows, the errors also grow. When enough samples are used, smaller
  $\delta$ leads to more accurate solutions. Finally, when not enough samples
  are used, larger $\delta$ helps to control the errors.}  
\label{fig:hj_highlights}
\end{figure}

Figure \ref{fig:hj_highlights} shows the median along with the 20th and 80th
percentiles of HJ residuals for the Hamiltonians with $p \in \{2,5,10\}$ in 
dimension $d=2$, over 50 repetitions. As expected, more samples lead to more
accurate solutions; when enough samples are used, the error vanishes as $\delta
\to 0^+$ (the green curve diminishes at the left end of each plot). Perhaps more
interesting is the effect of $\delta$ when \emph{not} enough samples are used:
in this case we find that \emph{larger} values of $\delta$ will often lead to
more accurate solutions (the blue and red curves \emph{decrease} as $\delta$ 
increases, for small values of $\delta$). This happens because a larger $\delta$
has a greater regularization effect, via the viscous HJ PDE interpretation
discussed in Section \ref{sec:viscous_burgers} (and is consistent with what is 
observed in previous work, such as \citet{chaudhari2018deep,
  osher2023hamilton}). 

Appendix \ref{app:hj_experiments} displays the full set of results over all
norms $p$ and dimensions $d$ considered. The behavior is broadly similar to what
is observed and described above, but unsurprisingly, higher dimensions lead to
larger errors (they would require more samples). Also, $p = 1.1$ acts as
somewhat of an exceptional case, as it tends to be more difficult overall. This
is probably due to instability in the autodifferentiation calculation used for
the residual \eqref{eq:lapprox_hj_res} of \smash{$u^{\delta,N}(x,t)$} in
\eqref{eq:lapprox_hj_sol}. Recall, here $g$ is based on the conjugate $H^* = 
\|\cdot\|_q^q / q$ of $H = \|\cdot\|_p^p / p$, and as $1/p + 1/q = 1$, we have
$q \to \infty$ as $p \to 1$. 

\subsection{Proximal point algorithm}
\label{sec:prox_point}

We consider the use of Laplace's approximation of the proximal mapping
within the context of running the proximal point algorithm on a set of benchmark
functions \citep{hansen2009real}, which in all cases except one (the ``sphere''  
function) do not admit analytic proximal maps. We note that the proximal  
point algorithm \eqref{eq:prox_point} with \eqref{eq:gaussian_samples}, 
\eqref{eq:laplace_prox_avg} in place of the exact proximal operator, repeats 
the following update for $k = 1,2,3,\dots$:  
\begin{gather}
\label{eq:lapprox_pp_samples}
Y_i \iidsim N(x_{k-1}, \delta \lambda I), \; i = 1,\dots, N, \\ 
\label{eq:lapprox_pp_update}
x_k = \frac{\sum_{i=1}^N Y_i \exp(-f(Y_i) / \delta)}
{\sum_{i=1}^N \exp(-f(Y_i) / \delta)}.  
\end{gather}
This is very intuitive: we explore the space locally by sampling points
\eqref{eq:lapprox_pp_samples} in a neighborhood of our current iterate
$x_{k-1}$, and then we take our next iterate $x_k$ to be a weighted average
\eqref{eq:lapprox_pp_update} of these sample, where we exponentially tilt in
favor of samples with smaller criterion values.
 
\subsubsection{Benchmark functions}

\begin{figure}[p]
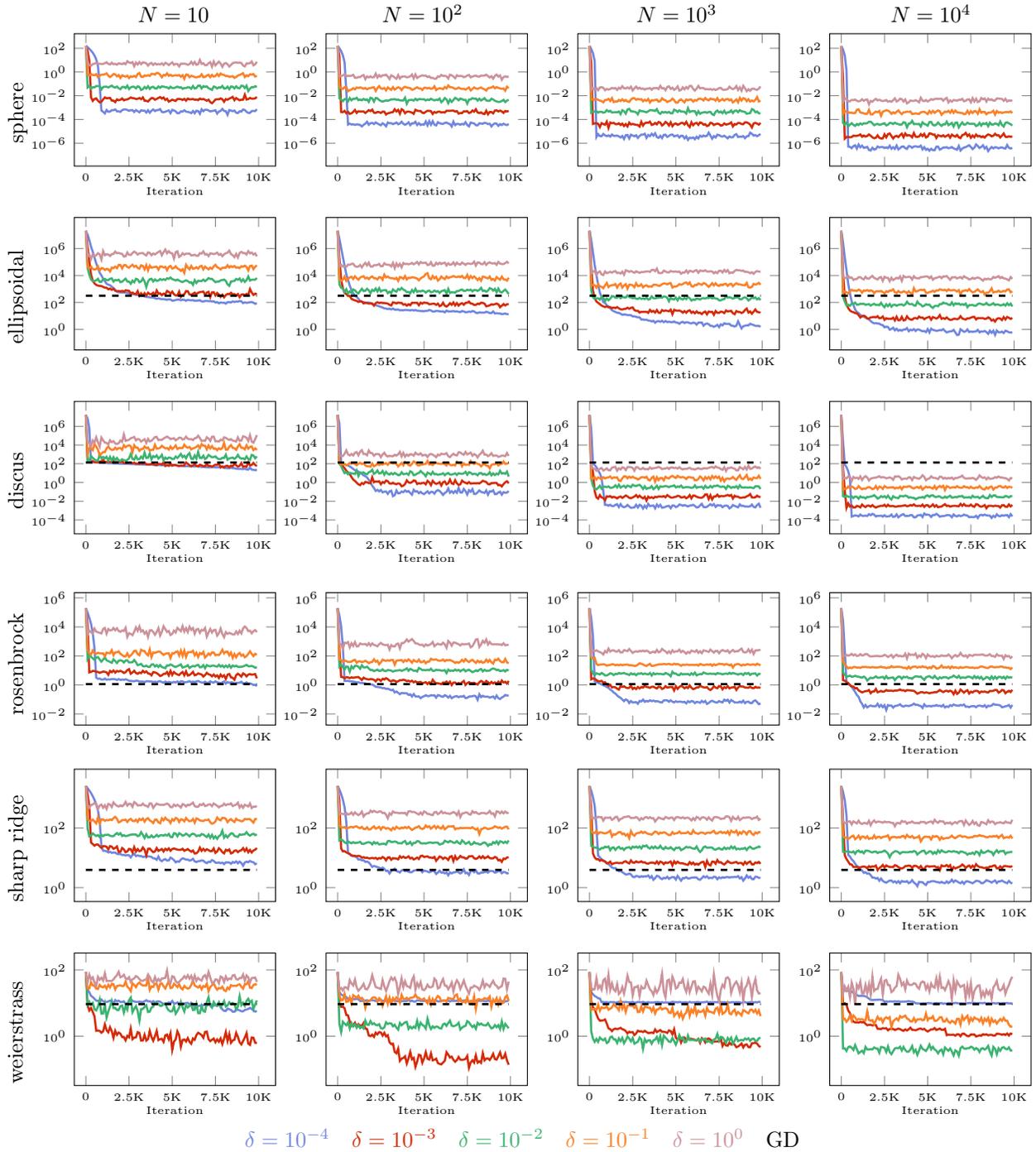

\input{images/prox_point_csv/prox-point-plots-list.tex}
\hspace*{0.82in}$N=10$
\hspace*{1.05in}$N=10^2$ 
\hspace*{1.05in}$N=10^3$ 
\hspace*{1.05in}$N=10^4$\\[-10pt]

\newcounter{ctrPlot}
\setcounter{ctrPlot}{0}   
\foreach [count=\k] \labelname/\filename in \plotInfo
{ \stepcounter{ctrPlot}
  \rottext{\labelname} 
  \includegraphics[page=\thectrPlot]{images/prox_point_csv/prox-point-plots.pdf}
  \hspace*{-7.5pt}
  \stepcounter{ctrPlot}
  \includegraphics[page=\thectrPlot]{images/prox_point_csv/prox-point-plots.pdf}
  \hspace*{-7.5pt}
  \stepcounter{ctrPlot}
  \includegraphics[page=\thectrPlot]{images/prox_point_csv/prox-point-plots.pdf}
  \hspace*{-7.5pt}
  \stepcounter{ctrPlot}
  \includegraphics[page=\thectrPlot]{images/prox_point_csv/prox-point-plots.pdf}\\
}
\begin{center}
  \vspace*{-0.4in}
  \stepcounter{ctrPlot}
  \includegraphics[page=\thectrPlot]{images/prox_point_csv/prox-point-plots.pdf}
\end{center}

\caption{\small
  Laplace-proximal point algorithm applied to various benchmark criterion
  functions. This algorithm repeats the iterations
  \eqref{eq:lapprox_pp_samples}, \eqref{eq:lapprox_pp_update} for a
  particular noise level $\delta$ and number of samples $N$. Each row shows a
  different benchmark function, each column shows a different number of samples
  \smash{$N \in \{10, 10^2, 10^3, 10^4\}$}; and each panel display results for
  noise levels \smash{$\delta \in \{10^{-4}, 10^{-3}, 10^{-2}, 10^{-1},
    1\}$}. The result from running gradient descent (GD, tuned over both the
  choice of step size and added noise level, as explained in the text) is shown
  as a dashed black line. In the first row, which is a convex and
  well-conditioned example, GD obtains the global minimum criterion value of
  zero, and the dashed line is not visible. In all others, the Laplace-proximal 
  point method is able to compete with and even outperform GD.}    
\label{fig:benchmark_runs}
\end{figure}

In Figure \ref{fig:benchmark_runs}, we study the effect of various choices of 
\smash{$\delta \in \{10^{-4}, 10^{-3}, 10^{-2}, 10^{-1}, 1\}$}, \smash{$N
  \in \{10, 10^2, 10^3, 10^4\}$} in running the Laplace-proximal point method 
\eqref{eq:lapprox_pp_samples}, \eqref{eq:lapprox_pp_update} over functions
in the benchmark suite. Throughout we fix $\lambda=1$.  As a reference, in
each setting, we plot the lowest criterion value achieved by running gradient
descent (GD) over 10,000 iterations, as a dashed black line. To be as favorable
possible toward GD, for each setting, we tune over the choice of step size used
by GD, as well as a variance level for Gaussian noise to add to the gradient at
each iteration, which includes zero noise (usual GD). We do this because adding
noise may help escape local minima in some of the nonconvex benchmark
functions. We report the the \emph{best} result over all step sizes and noise
levels for GD (the one with the lowest criterion value 10,000 iterations) as the
dashed black line. All functions in the benchmark are in $d=10$ dimensions, and
obtain a global minimum criterion value of zero at the origin, $x^\star=0$. We 
initialize all algorithms at $x_0 = (4,\dots,4) \in \R^{10}$, and average all
results over 3 repetitions (of random noise generation used in the algorithms).

In the first setting  ``sphere'', the dashed black line is not visible, as
GD converges to the global minimum of zero, lying outside of the plotting
range. The ``sphere'' benchmark is the only one in which the criterion is both
convex and well-conditioned. Beyond the ``sphere'' example, we see that
Laplace-proximal point (LPP) can compete with and even clearly outperform GD. In
two settings, ``ellipsoidal'' and ``discus'', LPP starts to outperform GD at a
reasonably small number of samples $N$ and reasonably large noise level
$\delta$; while in two others, ``rosenbrock'' and``sharp ridge'', it only
outperforms GD for larger $N$ and smaller $\delta$. The ``weierstrass'' example,
meanwhile, is different: in contrast to all of the other benchmarks, LPP does
not show consistent improvement as $N$ increases and $\delta$ decreases.
This is likely due to the fact that this function has a complex landscape with
many local minima.   
 
\subsubsection{Comparison with RGF}

We now compare to a well-known zeroth-order method based on Gaussian smoothing,
proposed by \citet{nesterov2017random}. This method is based on gradient
descent, but approximates the gradient of the criterion $f$ at each iterate
$x_{k-1}$ using a finite difference approximation that employs Gaussian
perturbations:      
\begin{gather}
\label{eq:rgf_samples}
Y_i \iidsim N(x_{k-1}, \delta I), \; i = 1,\dots, N, \\ 
\label{eq:rgf_update}
x_k = x_{k-1} - \frac{\eta}{N} \sum_{i=1}^N \frac{f(Y_i) - f(x_{k-1})}{\delta},
\end{gather}
where $\eta>0$ is a step size. \citet{nesterov2017random} refer to this
algorithm as the \emph{random gradient-free oracle} or RGF. (These authors only
consider $N=1$, but we find that averaging over multiple Gaussian draws tends to
improve results).       

While both use Gaussian sampling, it is interesting to note that the
motivation for RGF is quite different from that for Laplace-proximal point (LPP)
in \eqref{eq:lapprox_pp_samples}, \eqref{eq:lapprox_pp_update}. LPP can be
interpreted as follows:
\begin{itemize}
\item we first smooth $f$ using its Moreau envelope $f_\lambda$ (recalling that
  proximal point \eqref{eq:prox_point} is nothing more that gradient descent on
  the Moreau envelope \eqref{eq:moreau_grad_desc});
\item we then numerically approximate a gradient step with respect to the Moreau
  envelope (proximal map) using Laplace's method and Gaussian sampling.
\end{itemize}
On the other hand, RGF (as \citet{nesterov2017random} show) can be interpreted
as follows: 
\begin{itemize}
\item we first smooth $f$ using a Gaussian convolution;
\item we then numerically approximate a gradient of this convolved function
  using Gaussian sampling.  
\end{itemize}

\begin{figure}[t]
\centering
\begin{tabular}{ccc}
\hspace*{-8pt}
\includegraphics[width=0.32\textwidth, page=1]{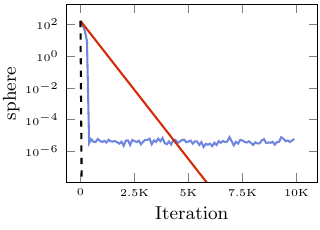} &
\hspace*{-8pt}
\includegraphics[width=0.32\textwidth, page=2]{images/rgf_comparisons/rgf-comparisons.pdf} &
\hspace*{-8pt}
\includegraphics[width=0.32\textwidth, page=3]{images/rgf_comparisons/rgf-comparisons.pdf} \\
\hspace*{-8pt}
\includegraphics[width=0.32\textwidth, page=4]{images/rgf_comparisons/rgf-comparisons.pdf} &
\hspace*{-8pt}
\includegraphics[width=0.32\textwidth, page=5]{images/rgf_comparisons/rgf-comparisons.pdf} &
\hspace*{-8pt}
\includegraphics[width=0.32\textwidth, page=6]{images/rgf_comparisons/rgf-comparisons.pdf} 
\end{tabular}
\caption{\small
  Comparison of \textcolor{TypalRed}{random-gradient free oracle (RGF) in red},
  as in \eqref{eq:rgf_samples}, \eqref{eq:rgf_update};
  \textcolor{TypalBlue}{Laplace-proximal point (LPP) in blue}, as in
  \eqref{eq:lapprox_pp_samples} \eqref{eq:lapprox_pp_update}; and gradient
  descent (GD) in dotted black; on the same set of benchmark critera from Figure   
  \ref{fig:benchmark_runs}. We fix $N=1000$ samples for RGF and LPP, fix
  $\lambda=1$ for LPP, and otherwise allow each algorithm to tune over their
  respective tuning parameters (noise level $\delta$ for each algorithm, and
  step size $\eta$ for RGF and GD).}    
\label{fig:rgf_comparison}
\end{figure}

Figure \ref{fig:rgf_comparison} compares RGF and LPP on the same set of
benchmark criteria as in Figure \ref{fig:benchmark_runs}. We fix
$N=1000$ for each algorithm, to equalize their sampling cost. Unlike Figure 
\ref{fig:benchmark_runs}, we now tune LPP over the noise level $\delta$,
reporting results for the best noise level (resulting in the smallest criterion
value in 10,000 iterations) in each setting. To be as favorable as possible to
RGF, we tune it over \emph{both} the noise level $\delta$ and the step size
$\eta$, and report results for the best combination of noise level and step size
in each benchmark. The results for GD are also shown, where we again tune it
over both the noise level and step size, as explained previously.  
Outside of the convex and well-conditioned ``sphere'' benchmark, where GD
performs best, RGF sometimes improves on GD (``sharp ridge'', ``weierstrass''), 
and other times RGF and GD perform quite similarly (``ellipsoidal'', ``discus'',
``rosenbrock''). LPP is competitive overall, sometimes improving on RGF
(``ellipsoidal'', ``discus'', ``rosenbrock'') and sometimes not (``sphere'',
``weierstrass'').        

\subsection{Bregman proximal gradient descent}
\label{sec:bregman_prox}

\begin{figure}[t]
\centering
\subfloat[Well-conditioned]{
\includegraphics[width=0.45\textwidth, page=1]{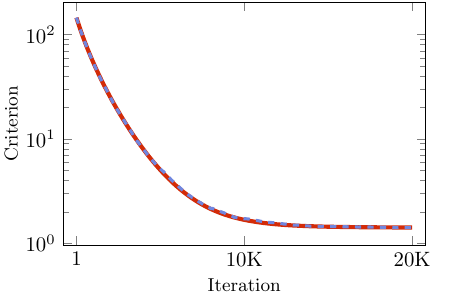}}
\subfloat[Ill-conditioned]{
\includegraphics[width=0.45\textwidth, page=2]{images/bpgd/bpgd-plots.pdf}}
\caption{\small
Comparison of BPGD approaches in Poisson linear inverse problems where the
underlying linear transform is well-conditioned or ill-conditioned. Here we
plot \textcolor{TypalRed}{exact BPGD with Burg entrope in red}, as in
\eqref{eq:bpgd2}; \textcolor{TypalBlue}{Laplace-BPGD with Burg entropy in dotted  
  blue}, as in \eqref{eq:bpgd1} where \smash{$h(x) = -\sum_{i=1}^d \log x_i$};
and \textcolor{TypalGreen}{Laplace-BPGD with the variable-metric majorizer in 
  green}, as in \eqref{eq:bpgd1} where \smash{$h(x) = -\sum_{i=1}^d \log 
  (Ax)_i$}. The Laplace methods use $\delta = 2 \cdot 10^{-3}$ as the noise
level and $N = 5 \cdot 10^4$ as the number of samples. We can see that the
Laplace-BPGD method with Burg majorizer tracks its exact counterpart closely. 
In the ill-conditioned setting, the variable-metric majorizer accelerates
convergence.}      
\label{fig:bpgd}
\end{figure}

Lastly, we consider the use of Laplace's method to approximate a Bregman
proximal map within the context of Bregman proximal gradient descent (BPGD). In
particular, we examine a variation on BPGD proposed by
\citet{bauschke2017descent}, which uses an elegant majorization scheme to handle
composite criterion functions where the differentiable part lacks a Lipschitz
continuous gradient. We focus on a regularized Poisson linear inverse problem,
as studied in their paper, where we aim to minimize  
\begin{equation}
\label{eq:poisson_inverse}
f(x) = D_\phi(b, Ax) + \mu \| x \|_1,
\end{equation}
over \smash{$\R^d_+$}, where \smash{$b \in \R^n_+$}, \smash{$A \in \R^{n \times
    d}_{++}$}, $\mu \geq 0$ is a regularization parameter, $\phi$ is the
Boltzmann-Shannon entropy \smash{$\phi(x) = \sum_{i=1}^d x_i \log x_i$}, and 
$D_\phi$ denotes its corresponding Bregman divergence,        
\begin{align}
\nonumber 
D_\phi(b, Ax) 
&= \phi(b) - \phi(Ax) - \nabla \phi(Ax)^\T (b-Ax) \\ 
\label{eq:poisson_bregman}
&= \underbrace{\sum_{i=1}^n\bigg( b_i \log \frac{b_i}{(Ax)_i} + b_i - (Ax)_i
  \bigg)}_{d(x)}.   
\end{align}
To minimize \eqref{eq:poisson_inverse}, we consider the BPGD algorithm of
\citet{bauschke2017descent} and choose as the majorizer the Burg entropy
\smash{$h(x) = -\sum_{i=1}^d \log x_i$}. For the Poisson inverse problem, the 
Burg entropy is a specially-crafted majorizer which satisfies two key
properties:    
\begin{itemize}
\item it majorizes the Bregman divergence in \eqref{eq:poisson_bregman}, in the
  sense that $Lh - d$ is convex for some $L>0$, which
  \citet{bauschke2017descent} prove is true for any choice $L \geq \|b\|_1$;    
\item it admits a closed-form Bregman proximal update (in the second line
  below), 
\begin{align}
\label{eq:bpgd1}
x_k 
&= \argmin_x \, \bigg\{ \mu \|x\|_1 + \nabla d(x_{k-1})^\T x + \frac{1}{\eta}
  D_h(x,x_{k-1}) \bigg\} \\  
\label{eq:bpgd2}
&= \bigg[ \frac{x_{k-1,i}}{1 + \eta (\mu + \nabla_i d(x_{k-1})) \cdot x_{k-1,i}}
  \bigg]_{i=1}^d,
\end{align}
were $\eta \in (0, 1/L)$ is a step size.
\end{itemize} 
The left panel of Figure \ref{fig:bpgd} shows an example of running the exact
BPGD update \eqref{eq:bpgd2}, and Laplace's method to approximate the general
form in \eqref{eq:bpgd1}, which does not use knowledge of the fact that for
the Burg entropy $h$ the Bregman proximal mapping is exact. For the Laplace   
approximation, we are able to carry out this out \emph{per coordinate}, because
the Bregman proximal mapping separates into a minimization problem per
coordinate, due to the separability of the Burg entropy itself. For each
coordinate, we then use importance sampling as in \eqref{eq:q_samples},
\eqref{eq:laplace_inf_prox_avg_iw}, with $q$ being the uniform density over  
$[10^{-6}, 50]$. 

In the left panel, we set $n=d=5$, and generate \smash{$A \in \R^{5 \times
    5}_{++}$} by sampling its entries independently from a uniform distribution
on $[1,2]$. We generate \smash{$\bar{x} \in \R^5$} by sampling its entries 
independently from a uniform on $[5,6]$, and randomly set half of these to 0.
Then, we generate \smash{$b \in \R^5_+$} by sampling its entries independently
from Poisson distributions with means \smash{$(A\bar{x})_i$}, $i=1,\dots,5$. We
fix $\mu=10^{-3}$ for the regularization parameter, $\eta=10^{-5}$ for the step
size, $\delta = 2 \cdot 10^{-3}$ for the level of noise, and $N = 5 \cdot 10^4$
for the number of samples. The figure shows the convergence curves (criterion
values per iteration) for the exact and Laplace-based BPGD methods. These look
identical, as should be expected for such a large number of samples $N$ and
small noise level $\delta$ in $d=5$ dimensions.

Meanwhile, the right panel of Figure \ref{fig:bpgd} shows a different example
setting in which $A$ is ill-conditioned and BPGD with the Burg entropy (whether
exact or Laplace-based) converges slowly as a result. The setup is the same 
except that we set $A = a a^\T$, where \smash{$a \in \R^5_{++}$} has entries
sampled independently from a uniform on $[0,1]$. Now, in addition to the exact
BPGD and Laplace-BPGD methods that use the Burg entropy, we consider a third
method: we use the \emph{variable-metric} majorizer \smash{$h(x) = -\sum_{i=1}^d 
  \log (Ax)_i$}. This does not lead to an exact proximal update, but nonetheless
\eqref{eq:bpgd1} can be approximated using Laplace's method and importance
sampling, as in \eqref{eq:q_samples}, \eqref{eq:laplace_inf_prox_avg_iw}, where 
we take $q$ to be the exponential distribution, arising from the term $\exp(-\mu
\|x\|_1)$ that appears in the integral. As we can see in the right panel of the
figure, such a variable-metric Laplace-BPGD approach converges more rapidly  
than standard BPGD, because the variable-metric approach has effectively
transformed the parameter space to account for the underlying geometry. 

\section{Discussion} 

In this work, we study connections between Laplace's approximation and infimal  
convolutions, with an eye toward solving optimization problems (using
proximal-type methods) as well as Hamilton-Jacobi PDEs. Our theory reflects the 
broad applicability of Laplace's method for approximating infimal convolutions,
allowing us to handle nonconvex functions, and settings where a minimizer lies
on the boundary of the domain (such as projections). Our experiments demonstrate
the versatility of the Laplace approximation in a few different problem areas,
which span optimization and PDEs. Practically, the challenge of sampling in
high-dimensional spaces remains a significant hurdle, and one that we do not
attempt to address at all. Any implementation based on Laplace's method which
strives for both precision and efficiency should likely invest in more
advanced sampling techniques. Future work may be able to provide samplers
tailored to particular environments and end-goals in optimization and PDEs.   
An open-source repository with code to reproduce our experiments is available at
\url{https://github.com/mines-opt-ml/laplace-inf-conv}. 

\subsection*{Acknowledgments}

We thank Lawrence Craig Evans and Luis Tenorio for helpful discussions.  This
work was supported in part by Office of Naval Research MURI grant
N00014-20-1-2787, to RJT and SO; Air Force Office of Scientific Research MURI
grant FA9550-18-502, to SO; and National Science Foundation grant DMS211074, to 
SWF. 

\bibliographystyle{plainnat}
\bibliography{ryantibs, lapprox}

\clearpage
\appendix

\section{Proofs}

\subsection{Proof of Theorem \ref{thm:laplace_method}}
\label{app:laplace_method}

For ease of reference, we state the theorem before its proof.

\medskip
\noindent
\textbf{Theorem \ref{thm:laplace_method}.}
\textit{\thmLaplaceMethod{app}}

\begin{proof}
Without loss of generality, assume that $x^\star = 0$, $\varphi(x^\star) = 0$,
and the Hessian of $\varphi$ at $x^\star$ is the identity matrix, \ie $\nabla^2
\varphi(x^\star) = I$. The first equality $x^\star = 0$ can be achieved by
shifting the domain, and the second $\varphi(x^\star) = 0$ can be achieved by
centering the function, which contributes the factor of
$\exp(\varphi(x^\star/\delta))$ to the left-hand side in
\eqref{eq:laplace_method_app}. The Hessian condition can be obtained via a
change of variables, and this contributes the factor \smash{$\sqrt{\det(\nabla^2
    \varphi(x^\star))}$} to the left-hand side in
\eqref{eq:laplace_method_app}. Now define 
\[
I_\delta = \frac{1}{(2\pi\delta)^{d/2}} \int_\cK h(x) \exp(-\varphi(x) / \delta)
\, \d{x}, \quad \text{for $\delta>0$}. 
\]
We seek to show that
\[
\lim_{\delta \to 0^+} I_\delta = h(0).
\]
We proceed by separately considering whether the set $\cK$ is compact or not. 
    
\paragraph{Compact case.} 

Suppose $\cK$ is compact. Following \citet{bach2021approximating}, define the
function $g\colon \cK \to \R$ by  
\[
g(x) = \begin{cases}
\begin{array}{cl}
\dfrac{\varphi(x)}{\|x\|_2^2} & \text{if $x \not= 0$}, \\[8pt]
1/2 & \text{otherwise}.
\end{array}
\end{cases}
\]
First note that $g$ is continuous at each $x \not= 0$ since $\varphi$ and 
$\|\cdot\|_2$ are both continuous (and the denominator is nonzero). Moreover, we
can show that $g$ is continuous at $x = 0$ via a second-order Taylor expansion
of $\varphi$:      
\[
\varphi(x) = \varphi(x^\star) + \nabla \varphi(x^\star)^\T x + \frac{1}{2} x^\T 
\nabla^2 \varphi(x^\star) x + R(x) = \frac{1}{2} \|x\|_2^2 + R(x), 
\]
where the remainder term is such that $R(x)/\|x\|_2^2 \to 0$ as $x \to 0$, and
therefore,      
\[
\lim_{x \to 0} g(x) = \lim_{x\to 0} \frac{\frac{1}{2} \|x\|_2^2 +
  R(x)}{\|x\|_2^2} = \frac{1}{2} + 0 = g(0). 
\]
Hence we have shown $g$ is continuous on all of $\cK$. 

Using the definition of $g$, the integral $I_\delta$ may be equivalently
written as  
\begin{align}
\nonumber
I_\delta 
&= \frac{1}{(2\pi\delta)^{d/2}} \int_\cK h(x) \cdot \exp\Big(-\|x\|_2^2 \cdot
  g(x) / \delta\Big) \, \d{x} \\   
\label{eq:I_delta_rewrite}
&= \frac{1}{(2\pi)^{d/2}} \int_{\R^d} \one_\cK(\sqrt\delta y) \cdot
  h(\sqrt\delta y) \cdot  \exp\Big(-\|y\|_2^2\cdot g(\sqrt\delta y)\Big) \,
  \d{y},  
\end{align}
where the second line follows from the change of variables \smash{$x =
  \sqrt\delta y$} and $\one_\cK$ is the indicator function of $\cK$,  
\[
\one_\cK(x) = 
\begin{cases}
1 & \text{if $x \in \cK$} \\  
0 & \text{otherwise}.
\end{cases}
\]
For each $y$, the absolute value of the integrand in \eqref{eq:I_delta_rewrite}
is upper bounded by    
\begin{equation}
\Big( \max_{x \in \cK} \, |h(x)| \Big) \cdot \exp\bigg(-\|y\|_2^2 \cdot \Big( 
\min_{x \in \cK} \, g(x) \Big) \bigg). 
\label{eq:dominating_term}
\end{equation}
By compactness of $\cK$ and continuity of $h$, the max term is finite. By
compactness of $\cK$, continuity of $g$, and the fact that $g$ is strictly
positive on $\cK$, the min term is strictly positive. These facts imply    
\[
\int_{\R^d} \Big( \max_{x \in \cK} \, |h(x)| \Big) \cdot \exp\bigg(-\|y\|_2^2 
\cdot \Big( \min_{x \in \cK} \, g(x) \Big) \bigg)\, \d{y} < \infty,
\]
so the integrand in \eqref{eq:I_delta_rewrite} is dominated by an integrable
function \eqref{eq:dominating_term}. By the dominated convergence theorem, 
\begin{align*}
\lim_{\delta \to 0^+} I_\delta 
&=  \frac{1}{(2\pi)^{d/2}} \int_{\R^d} \lim_{\delta \to 0^+} \one_\cK 
(\sqrt\delta y) \cdot h(\sqrt\delta y) \cdot \exp\Big(-\|y\|_2^2 \cdot
g(\sqrt\delta y)\Big) \, \d{y} \\ 
&= h(0) \cdot \frac{1}{(2\pi)^{d/2}} \int_{\R^d} \exp(-\|y\|_2^2/2) \, \d{y} \\  
&= h(0),
\end{align*}
completing the proof of the compact case.
    
\paragraph{Noncompact case.}

Suppose $\cK$ is not compact. In this case, we decompose the integral $I_\delta$
into two parts:     
\[
I_\delta = \underbrace{\frac{1}{(2\pi\delta)^{d/2}} \int_{\cS_\epsilon} h(x)
  \exp(-\varphi(x) / \delta) \, \d{x}}_{I_{1,\delta}} \,+\,
\underbrace{\frac{1}{(2\pi\delta)^{d/2}} \int_{\cK \setminus \cS_\epsilon} h(x)
  \exp(-\varphi(x) / \delta) \, \d{x}}_{I_{2,\delta}}.  
\]
Since $\varphi$ is continuous, the set $\cS_\epsilon$ is closed. Thus, because
it is also bounded (by assumption), the set $\cS_\epsilon$ is compact and the
arguments from the last case show \smash{$I_{1,\delta} \to h(0)$} as $\delta \to 
0$. To complete the proof, it remains to show \smash{$I_{2,\delta} \to 0$} as
$\delta \to 0^+$. Bringing the leading term inside the integral and taking 
absolute values reveals:
\[
|I_{2,\delta}| \leq \int_{\R^d} \underbrace{\one_{\cK \setminus \cS_\epsilon}(x)
  \cdot \frac{|h(x)|}{(2\pi\delta)^{d/2}} \cdot \exp(-\varphi(x) /
  \delta)}_{f_\delta(x)} \, \d{x}. 
\]
Furthermore, for $x \in \cK \setminus \cS_\epsilon$,  
\begin{align*}
\frac{\d{f}_\delta(x)}{\d\delta}
&= \frac{|h(x)|}{(2\pi \delta)^{d/2}} \cdot \bigg( \frac{\varphi(x)}{\delta^2}
  - \frac{d}{2\delta} \bigg) \cdot \exp(-\varphi(x)/\delta) \\
&\geq \frac{|h(x)|}{(2\pi \delta)^{d/2}} \cdot \bigg( \frac{\epsilon}{\delta^2}
  - \frac{d}{2\delta} \bigg) \cdot \exp(-\varphi(x)/\delta).
\end{align*} 
If $\delta \leq 2\epsilon/d$, then the lower bound in the last line above is 
nonnegative; in other words, we have shown that for $\delta \leq 2\epsilon/d$,
the quantity $f_\delta(x)$ is nonincreasing as $\delta$ decreases, at each $x
\in \cK \setminus \cS_\epsilon$. As $f_\delta \geq 0$, we conclude that
$f_\delta$ is dominated by $f_{2\epsilon/d}$, which is integrable by
assumption, and another application of the dominated convergence theorem
therefore yields     
\begin{align*}
\lim_{\delta \to 0^+} |I_{2,\delta}|
&\leq \int_{\R^d} \lim_{\delta \to 0^+} f_\delta(x) \, \d{x} \\
&\leq \int_{\cK \setminus \cS_\epsilon} \lim_{\delta \to 0^+} 
\frac{|h(x)|}{(2\pi \delta)^{d/2}} \cdot \exp(-\epsilon/\delta) \, \d{x} \\
&=0,  
\end{align*}
where the last line holds since $\delta^{-d/2} \exp(-\epsilon / \delta) \to 0$
as $\delta \to 0^+$. This completes the proof of the noncompact case, and the 
theorem.     
\end{proof}

\subsection{Proof of Theorem \ref{thm:laplace_sn}}
\label{app:laplace_sn}

For ease of reference, we state the theorem before its proof.

\medskip
\noindent
\textbf{Theorem \ref{thm:laplace_sn}.}
\textit{\thmLaplaceSelfNormalized{app}}

\begin{proof}
As before, we assume without loss of generality that $x^\star = 0$ and
$\varphi(x^\star) = 0$. We also assume that the local H{\"o}lder bound in
\eqref{eq:local_holder_app} holds on all of $\cS_\epsilon$; this is possible
because $\{\cS_t : t \leq \epsilon\}$ is a family of compact sets which are
decreasing according to the partial ordering given by set inclusion, and $\cS_0
= \{x^\star\}$, thus if needed we can just redefine $\epsilon$ to be small
enough such that $\cS_\epsilon \subseteq \cK \cap U$. Now let          
\[
  p^\delta(x) = \frac{\exp(-\varphi(x) / \delta)}{\int_\cK \exp(-\varphi(x) /  
  \delta) \, \d{x}}.
\]
For $\delta \leq \epsilon$, note that the denominator here is finite because
the integrand is upper bounded by $\exp(-\varphi(x) / \epsilon)$, which is 
integrable by assumption. Furthermore, the denominator is positive because the  
integrand is lower bounded by $\exp(-a \|x-x^\star\|_2^q) > 0$ on
$\cS_\epsilon$. Hence we have shown \smash{$0 < \int_\cK \exp(-\varphi(x) /
  \delta) \, \d{x} < \infty$}, which means that \smash{$p^\delta$} is a
well-defined probability measure on $\cK$. We now seek to prove
\[
\lim_{\delta \to 0^+} \int_\cK h(x) p^\delta(x) \, \d{x} = h(0),
\]
which holds if  
\[
\lim_{\delta \to 0^+} \underbrace{\int_\cK |h(x) - h(0)| p^\delta(x) \,
  \d{x}}_{J_\delta} = 0, 
\]
We split the integral $J_\delta$ into two parts:
\[
  J_\delta = \underbrace{\int_{\cK \cap B(0,\tau)} |h(x) - h(0)| p^\delta(x) \, 
  \d{x}}_{J_{1,\delta,\tau}} \,+\, \underbrace{\int_{\cK \setminus B(0,\tau)} 
  |h(x) - h(0)| p^\delta(x) \, \d{x}}_{J_{2,\delta,\tau}},
\]
where $B(0,\tau) = \{ x \in \R^d : \|x\|_2 \leq \tau \}$ is the ball of radius
$\tau>0$ centered at the origin, and the radius $\tau$ is yet to be specified.     
Let $\mu>0$ be given. It suffices to verify the existence of $\tau>0$ and
$\delta_0>0$ such that  
\[
J_{1,\delta,\tau} \leq \frac{\mu}{2} \quad \text{and}\quad 
J_{2,\delta,\tau} \leq \frac{\mu}{2}, \quad \text{for $\delta \leq \delta_0$}.   
\]
This would lead to \smash{$J_\delta = J_{1,\delta,\tau} + J_{2,\delta,\tau} \leq
  \mu$} for all $\delta \leq \delta_0$, implying the desired convergence
$J_\delta \to 0$ as $\delta \to 0^+$.  The rest of this proof is structured as
follows. We bound \smash{$J_{1,\delta,\tau} \leq \mu/2$} by choosing $\tau$ to
be sufficiently small (step 1). This bound holds for any $\delta>0$. For this
same $\tau$, we upper bound \smash{$J_{2,\delta,\tau}$} by a quantity that 
depends on $\delta$ (step 2). We then show that the numerator of this upper
bound converges to zero (step 3), whereas the denominator is bounded below by a
positive constant as $\delta \to 0^+$ (step 4). Together, steps 3 and 4 imply
the existence of $\delta_0$ such that \smash{$J_{2,\delta,\tau}\leq \mu/2$} for 
all $\delta \leq \delta_0$, which will complete the proof.

\paragraph{Step 1: bounding the integral $J_{1,\delta}$.} 

As $h$ is continuous, we can choose $\tau$ so that $|h(x) - h(0)| \leq \mu/2$
for all $x \in B(0,\tau)$, which makes   
\[
J_{1,\delta,\tau} \leq \frac{\mu}{2} \cdot \int_{\cK \cap B(0,\tau)}
p^\delta(x) \, \d{x} \leq \dfrac{\mu}{2}. 
\]

\paragraph{Step 2: bounding the integral $J_{2,\delta}$.} 

Observe
\begin{align} 
\nonumber
 J_{2,\delta,\tau} 
&= \frac{\delta^{-d/q} \cdot \int_{\cK \setminus B(0,\tau)} |h(x) - h(0)|
  \exp(-\varphi(x) / \delta) \, \d{x}}{\delta^{-d/q} \cdot \int_\cK 
  \exp(-\varphi(x) / \delta) \, \d{x}} \\
\label{eq:j2_delta_bound}
& \leq \frac{\delta^{-d/q} \cdot \int_{\cK \setminus B(0,\tau)} |h(x) - h(0)|
  \exp(-\varphi(x) / \delta) \, \d{x}}{\delta^{-d/q} \cdot \int_{\cK \cap
  B(0,\tau)} \exp(-\varphi(x) / \delta) \, \d{x}}. 
\end{align} 
The motivation for introducing the factor of $\delta^{-d/q}$ in the numerator
and denominator above is a change of variables that will be used in step 4.  

\paragraph{Step 3: analyzing the numerator in \eqref{eq:j2_delta_bound}.} 

Let $\epsilon_0>0$ be such that \smash{$\cS_{\epsilon_0} \subseteq
  \cK \cap B(0,\tau)$}. Such a value of $\epsilon_0$ always exists because,
similar to an argument given earlier, $\{\cS_t : t \leq \epsilon \}$ is a family
of nested compact sets with $\cS_0 = \{x^\star\}$. We assume without loss of
generality that $\epsilon_0 \leq \epsilon$ (otherwise we just make $\epsilon_0$
smaller). Now rewrite the numerator in \eqref{eq:j2_delta_bound} as          
\[
\delta^{-d/q} \cdot \int_{\cK \setminus B(0,\tau)} |h(x) - h(0)| 
  \exp(-\varphi(x) / \delta) \, \d{x} = \int_{\R^d} \underbrace{\one_{\cK
  \setminus B(0,\tau)}(x) \cdot \frac{|h(x) - h(0)|}{\delta^{d/q}} \cdot
  \exp(-\varphi(x) / \delta)}_{f_\delta(x)} \, \d{x}.
\]
Similar to a calculation given at the end of the proof of Theorem
\ref{thm:laplace_method} in Appendix \ref{app:laplace_method}, for $x \in \cK
\setminus B(0,\tau)$,     
\begin{align*}
\frac{\d{f}_\delta(x)}{\d\delta}
&= \frac{|h(x) - h(0)|}{\delta^{d/q}} \cdot \bigg( \frac{\varphi(x)}
  {\delta^2} - \frac{d}{q\delta} \bigg) \cdot \exp(-\varphi(x)/\delta) \\  
&\geq \frac{|h(x) - h(0)|}{\delta^{d/q}} \cdot \bigg( \frac{\epsilon_0} 
  {\delta^2}- \frac{d}{q\delta} \bigg) \cdot \exp(-\varphi(x)/\delta), 
\end{align*}
where the last line holds since $\varphi(x) \geq \epsilon_0$ (which is implied
by our choice of $\epsilon_0$ such that \smash{$\cS_{\epsilon_0} \subseteq
  B(0,\tau)$}). Thus for all $\delta \leq q\epsilon_0/d$, the quantity
$f_\delta(x)$ is nonincreasing as $\delta$ decreases, at each $x \in  \cK
\setminus B(0,\tau)$. As $f_\delta \geq 0$, and $f_{q\epsilon_0/d}$ is
integrable (following from $\epsilon_0 \leq \epsilon$, and the fact that
$f_{q\epsilon/d}$ is integrable by assumption), we can apply the dominated
convergence theorem to yield     
\begin{align*}
\lim_{\delta \to 0^+} \int_{\R^d} f_\delta(x) \, \d{x}
&= \int_{\R^d} \lim_{\delta \to 0^+} f_\delta(x) \, \d{x} \\
&\leq \int_{\cK \setminus B(0,\tau)} \lim_{\delta \to 0^+} \frac{|h(x) -
   h(0)|}{\delta^{d/q}} \cdot \exp(-\epsilon_0/\delta) \, \d{x} \\  
&= 0,  
\end{align*}
where the last line holds since $\delta^{-d/q} \exp(-\epsilon_0 / \delta) \to 0$ 
as $\delta \to 0^+$. 
 
\paragraph{Step 4: analyzing the denominator in \eqref{eq:j2_delta_bound}.} 

Recall that $x^\star = 0$, and we that assume the minimizer is in the interior
of $\cK$, thus we may assume without loss of generality that $B(0,\tau)
\subseteq \cK$ (otherwise in step 1 we simply take $\tau$ to be smaller). We may
further assume that the local H{\"o}lder condition in
\eqref{eq:local_holder_app} holds on $B(0,\tau)$ (again, otherwise in step 1 we
just take $\tau$ to be smaller). We can then lower bound the denominator in
\eqref{eq:j2_delta_bound} as follows:
\begin{align*}
\delta^{-d/q} \cdot \int_{\cK \cap B(0,\tau)} \exp(-\varphi(x) / \delta) \,
  \d{x} 
&\geq \delta^{-d/q} \cdot \int_{B(0,\tau)} \exp(-a\|x\|_2^q/\delta) \, \d{x} \\  
&= \int_{B(0, \tau \delta^{-1/q})} \exp(-a\|y\|_2^q) \, \d{y} \\  
&= \int_0^{\tau \delta^{-1/q}} \int_{\partial B(0,r)} \exp(-ar^q) \, \d\cH^{d-1} 
  \, \d{r}, 
\end{align*}
where the second line follows from a change of variables \smash{$x =
  \delta^{1/q} y$}, and the last line from a change to polar coordinates (see,
\eg Appendix C.3 in \citet{evans2010partial}), with $\cH^{d-1}$ denoting
Hausdorff measure of dimension $d-1$. Now, $\cH^{d-1}(\partial B(0,r)) = c
r^{d-1}$ for a constant $c>0$ depending only on $d$, so    
\[
\int_0^{\tau \delta^{-1/q}} \int_{\partial B(0,r)} \exp(-ar^q) \, \d\cH^{d-1} \,
\d{r} = c \cdot \int_0^{\tau \delta^{-1/q}} r^{d-1} \exp(-ar^q) \, \d{r},    
\]
and therefore in the limit the denominator is lower bounded by 
\[
\lim_{\delta \to 0^+} c \cdot \int_0^{\tau \delta^{-1/q}} r^{d-1} \exp(-ar^q) \,
\d{r} = c \cdot \int_0^\infty r^{d-1} \exp(-ar^q) \, \d{r},  
\]
with the quantity in the last line being a positive constant. This completes the
proof of the theorem.     
\end{proof}

\subsection{Local H{\"o}lder condition for weakly convex functions} 
\label{app:local_holder_weakly_convex}

Recall, a function $\varphi$ is said to be $\rho$-weakly convex on a set $U$
provided that the map \smash{$x \mapsto \varphi(x) + (\rho/2) \|x\|_2^2$} is
convex on $U$. We now verify that if $f$ and $g$ are $\rho_f$-weakly convex and
$\rho_g$-weakly convex, respectively, on a neighborhood $U$ of $y_x$, then the
local H{\"o}lder condition \eqref{eq:local_holder} in Corollary
\ref{cor:laplace_inf_conv} holds for $\varphi_x$ with exponent $q=1$. Note first
that convexity of $f(y) + (\rho_f/2) \|y\|_2^2$ implies convexity of    
\[
f(y) - f(y_x) + \frac{\rho_f}{2} \|y\|_2^2 - 2 y^\T y_x + \frac{\rho_f}{2}
\|y_x\|_2^2 = \underbrace{f(y) - f(y_x) + \frac{\rho_f}{2} \|y -
  y_x\|_2^2}_{F_x(y)},   
\]
because we have only added a linear function and a constant. Similarly,
convexity of $g(y) + (\rho_g/2) \|y\|_2^2$ implies convexity of   
\[
g(x-y) - g(x-y_x) + \frac{\rho_g}{2} \|x-y\|_2^2 - 2 (x-y)^\T (x-y_x) +
\frac{\rho_g}{2}  \|x-y_x\|_2^2 = \underbrace{g(x-y) - g(x-y_x) +
  \frac{\rho_g}{2} \|y - y_x\|_2^2}_{G_x(y)},
\]
because we have only made an affine variable transformation $y \mapsto x-y$,
then added a linear function and a constant. Now, by definition 
\begin{align*}
\varphi_x(y) - \varphi_x(y_x) &= f(y) - f(y_x) + g(x-y) - g(x-y_x) \\
&= F_x(y) + G_x(y) - (\rho_f + \rho_g) \|y - y_x\|_2^2.
\end{align*}
As $F_x(y) + G_x(y)$ is convex, and convex functions are locally Lipschitz
(see, \eg Theorem 6.7 in \citet{evans2015measure}), we know that it is Lipschitz
on $U$. Furthermore, assuming without a loss of generality that $\diam(U) =
\sup\{ \|y - z\|_2 : y,z \in U \} \leq 1$ (otherwise we just shrink $U$ around
$y_x$ such that this holds), we have 
\[
(\rho_f + \rho_g) \|y - y_x\|_2^2 \leq (\rho_f + \rho_g) \|y - y_x\|_2,
\]
Thus $\varphi_x(y) - \varphi_x(y_x)$ is the sum of two Lipschitz functions on
$U$, and therefore it is itself Lipschitz on $U$, \ie locally H{\"o}lder with
exponent $q=1$. 

\subsection{Proof of Theorem \ref{thm:laplace_sn_bdry}}
\label{app:laplace_sn_bdry}

For ease of reference, we state the theorem before its proof.

\medskip
\noindent
\textbf{Theorem \ref{thm:laplace_sn_bdry}.}
\textit{\thmLaplaceSelfNormalizedBoundary{app}}

\begin{proof}
The proof is the same as that for Theorem \ref{thm:laplace_sn} in Appendix
\ref{app:laplace_sn}, except for step 4. Because $x^\star = 0$ is assumed to lie
on the boundary of $\cK$, it is no longer possible to take $\tau$ to be small
enough so that $B(0,\tau) \subseteq \cK$. However, we can still take $\tau$ to 
be small enough so that the local H{\"o}lder bound holds on $\cK \cap
B(0,\tau)$, which results in a lower bound for the denominator in
\eqref{eq:j2_delta_bound} of  
\begin{align}
\nonumber
\delta^{-d/q} \cdot \int_{\cK \cap B(0,\tau)} \exp(-\varphi(x) / \delta) \,
  \d{x}  
&\geq \delta^{-d/q} \cdot \int_{\cK \cap B(0,\tau)} \exp(-a\|x\|_2^q/\delta) \,
  \d{x} \\
\label{eq:denom_lower_bound}
&= \int_{\cK_\delta \cap B(0, \tau \delta^{-1/q})} \exp(-a\|y\|_2^q) \, \d{y},  
\end{align}
the second line using a change of variables \smash{$x = \delta^{1/q} y$}, where
we abbreviate \smash{$\cK_\delta = \delta^{-1/q} \cK_\delta = \{ \delta^{-1/q} x
  : x \in \cK \}$}. The rest of this proof proceeds as follows. We first show
that we can use polar coordinates to compute \eqref{eq:denom_lower_bound},   
despite the fact that the integral is over the (possibly) nonspherical set
\smash{$\cK_\delta \cap B(0, \tau \delta^{-1/q})$} (step 4a). We then show how
to lower bound this integral by a positive constant as $\delta \to 0^+$, using 
our assumptions on the set $\cK$ (step 4b). This will complete the proof.      

\paragraph{Step 4a: rewriting \eqref{eq:denom_lower_bound} using polar
  coordinates.} 

Rewrite \eqref{eq:denom_lower_bound} as 
\[
\int_{\cK_\delta \cap B(0, \tau \delta^{-1/q})}\exp(-a\|y\|_2^q) \, \d{y} =
\int_{B(0, \tau \delta^{-1/q})} \underbrace{\one_{\cK_\delta}(y) \cdot
  \exp(-a\|y\|_2^q)}_{\psi(y)}  \, \d{y}.  
\]
Let \smash{$v_k = \eta_{1/k} * \one_{\cK_\delta}$} be a mollified version of the
indicator function \smash{$\one_{\cK_\delta}$}, where $\eta_\epsilon$ is the
``standard'' mollifier with bandwidth $\epsilon>0$ (see, \eg Appendix C.5 in   
\citet{evans2010partial}). Then defining \smash{$\psi_k(y) = v_k(y) \cdot
\exp(-a\|y\|_2^q)$}, since mollified functions converge locally in $L^1$ (\eg
Theorem 7 in Appendix C.5 of \citet{evans2010partial}),  
\begin{align}
\nonumber
\int_{B(0, \tau \delta^{-1/q})} \psi(y) 
&= \lim_{k \to \infty} \int_{B(0, \tau \delta^{-1/q})} \psi_k(y) \\
\nonumber
&= \lim_{k \to \infty} \int_0^{\tau \delta^{-1/q}} \int_{\partial B(0,r)} v_k \exp(-ar^q) \,
  \d\cH^{d-1} \, \d{r} \\ 
\nonumber
&= \int_0^{\tau \delta^{-1/q}} \int_{\partial B(0,r)} \lim_{k \to \infty} v_k \exp(-ar^q) \,
  \d\cH^{d-1} \, \d{r} \\ 
\nonumber
&= \int_0^{\tau \delta^{-1/q}} \int_{\partial B(0,r)} \one_{\cK_\delta} \cdot
  \exp(-ar^q) \, \d\cH^{d-1} \, \d{r} \\ 
\label{eq:denom_lb_polar1}
&= \int_0^{\tau \delta^{-1/q}} \cH^{d-1}(\cK_\delta \cap \partial B(0,r)) \cdot
  \exp(-ar^q) \, \d{r}.      
\end{align}
The second line uses polar coordinates integration, which applies since $\psi_k$ 
is smooth (infinitely differentiable) by construction, the third uses the
dominated convergence theorem, which applies because $0 \leq v_k \leq 1$, and
the fourth uses the property that mollified functions converge pointwise almost  
everywhere (\eg Theorem 7 in Appendix C.5 of \citet{evans2010partial}). 

\paragraph{Step 4b: lower bounding \eqref{eq:denom_lb_polar1} as $\delta \to 
  0^+$.} 
\allowdisplaybreaks

Observe that the local star-shaped assumption \eqref{eq:local_star_app} can be   
reformulated equivalently (recalling that we have taken $x^\star = 0$) as  
\begin{equation}
\label{eq:local_star_reform} 
\beta \cK \cap \partial B(0,r) \supseteq \cK \cap \partial B(0,r), \quad
\text{for $r \leq r_0$ and $\beta \geq 1$}. 
\end{equation}
To see this, take $x \in \cK$ such that $\|x\|_2 = r \leq r_0$, and note that 
\eqref{eq:local_star_app} says $\alpha x \in \cK$ for $\alpha \leq 1$, or
equivalently, $x \in \beta \cK$ for $\beta = 1/\alpha \geq 1$. By taking $\delta 
\leq \min\{(\tau/r_0)^q, 1\}$, we have \smash{$\tau \delta^{-1/q} \geq r_0$}, so
we may lower bound the integral in \eqref{eq:denom_lb_polar1} by
\[
\int_0^{\tau \delta^{-1/q}} \cH^{d-1}(\cK_\delta \cap \partial B(0,r)) \cdot
  \exp(-ar^q) \, \d{r} \geq \int_0^{r_0} \cH^{d-1}(\cK_\delta \cap \partial
  B(0,r)) \cdot \exp(-ar^q) \, \d{r}.
\]
Applying \eqref{eq:local_star_reform} to the integrand (with \smash{$\beta =
  \delta^{-1/q} \geq 1$}) gives   
\begin{equation}
\label{eq:denom_lb_polar2}
\int_0^{r_0} \cH^{d-1}(\cK_\delta \cap \partial B(0,r)) \cdot \exp(-ar^q) \,
\d{r} \geq \int_0^{r_0} \cH^{d-1}(\cK \cap \partial B(0,r)) \cdot \exp(-ar^q) \,
\d{r}. 
\end{equation}
Meanwhile, for $r \leq r_0$, 
\begin{align}
\nonumber
\cH^{d-1}(\cK \cap \partial B(0,r)) 
&= \bigg( \frac{r}{r_0} \bigg)^{d-1} \cH^{d-1}\Bigg( \bigg( \frac{r_0}{r} \bigg)
  \cK \cap \partial B(0,r_0) \Bigg) \\
\label{eq:hausdorff_lb}
&\geq \bigg( \frac{r}{r_0} \bigg)^{d-1} \cH^{d-1}(\cK \cap \partial B(0,r_0)),
\end{align}
the first line using the $(d-1)$-homogeneity of $\cH^{d-1}$ (\eg Theorem 2.2 in
Chapter 2.1 of \citet{evans2015measure}), and the second line again using
\eqref{eq:local_star_reform}. 

We remark our assumption $\cL^d(\cK \cap \partial B(0,r_0)) > 0$, where $\cL^d$
denotes Lebesgue measure of dimension $d$, implies $\cH^{d-1}(\cK \cap \partial
B(0,r_0')) > 0$, for some $r_0' \leq r_0$. Otherwise, $\cH^{d-1}(\cK \cap
\partial B(0,r)) = 0$ for all $r \leq r_0$ would imply that $\cL^d(\cK \cap
\partial B(0,r_0)) = 0$, as we can represent the Lebesgue measure as an integral
over Hausdorff measure of boundary segments, by the same mollification argument
used to derive \eqref{eq:denom_lb_polar1} in step 4a.  Assume without a loss of
generality that $r_0' = r_0$ (otherwise we simply redefine $r_0$ to make this
true), and  abbreviate the lower bound in \eqref{eq:hausdorff_lb} by $c
r^{d-1}$, where $c>0$ is a constant depending only on $r_0$ and $d$. Then
applying \eqref{eq:hausdorff_lb} to \eqref{eq:denom_lb_polar2} gives 
\[
\int_0^{r_0} \cH^{d-1}(\cK \cap \partial B(0,r)) \cdot \exp(-ar^q) \, \d{r} \geq
c \cdot \int_0^{r_0} r^{d-1} \exp(-ar^q) \, \d{r},   
\]
the quantity in the last line being another positive constant. This completes 
the proof of the theorem.      
\end{proof}

\subsection{Local full-dimensional and star-shaped conditions for convex 
  bodies}  
\label{app:local_star_convex}

Recall, a convex body $\cK \subseteq \R^d$ is a closed convex set with nonempty
interior. Fix any $x^\star$ on the boundary of $\cK$. For any $x \in \cK$, we
have $\alpha x + (1-\alpha) x^\star \in \cK$ for all $\alpha \in [0,1]$ by
convexity, so clearly the local star-shaped condition \eqref{eq:local_star_main}
is met for any $r_0>0$.

Meanwhile, for any $r_0>0$, the set $\cK \cap B(x^\star,r_0)$ has nonempty
interior. To see this, take $x \in \interior(\cK)$ such that $\|x-x^\star\|_2
\leq r_0/3$. (This can be accomplished by taking an arbitrary point in the
interior, then shrinking toward $x^\star$, and invoking convexity, until the
distance bound is met.) Then by definition, there exists $\epsilon>0$ such that 
$B(x,\epsilon) \subseteq \interior(\cK)$. Defining $\epsilon' = \min\{\epsilon,
r_0/3\}$, we have that $B(x,\epsilon')$ is contained in
\[
\interior(\cK) \cap \interior(B(x^\star,r_0)) = \interior(\cK \cap
B(x^\star,r_0)),
\]
so $x$ lies in the interior of $\cK \cap B(x^\star,r_0)$. A nonempty interior
implies that $\cK \cap B(x^\star,r_0)$ has positive Lebesgue measure. This
verifies both conditions of Theorem \ref{thm:laplace_sn_bdry} for convex bodies.  
 
\clearpage
\section{Further Hamilton-Jacobi Experiments} 
\label{app:hj_experiments}

\begin{figure}[h] 
  \centering   
  \foreach [count=\m] \norm in {1.1,2,5,10}   
  { \hspace*{-12pt} 
    \rottext{\qquad\small Residual} 
    \hspace*{-18pt}
    \foreach [count=\k] \dim in {2,5,10}
    { \pgfmathtruncatemacro{\pageNum}{3 * (\m-1) + \k}
      \subfloat[Dimension $d=\dim$, Norm $p=\norm$]{
        \includegraphics[width=0.33\textwidth, page=\pageNum]{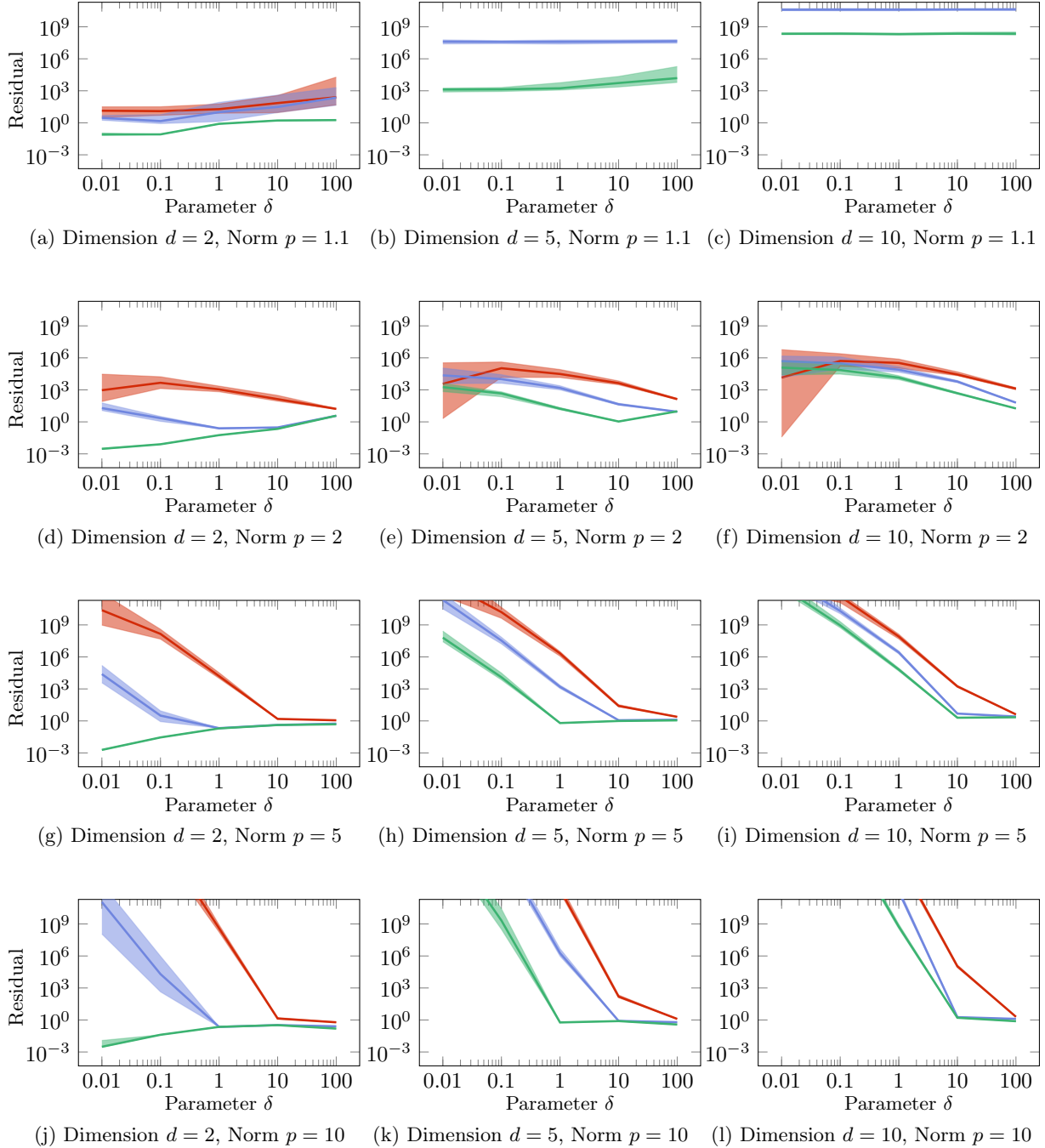}}
      \hspace*{-20pt}
      \ifthenelse{\k = 3}
      {\\[20pt]}
      {} 
    } 
  } 
\vspace{-10pt}
\caption{\small
  Residuals for HJ approximation in all dimensions $d$ and norms $p$ (for the
  Hamiltonian $H = \|\cdot\|_p^p / p$) that we consider. Each panel shows 
  \textcolor{TypalRed}{$N = 10$ samples via red}; \textcolor{TypalBlue}{$N=10^3$
    samples via blue}; and \textcolor{TypalGreen}{$N=10^5$ samples via green}. 
  In general, we see similar trends to what is observed in Figure
  \ref{fig:hj_highlights}, and additionally, we now see that as the dimension
  grows, the errors also grow. The case $p=1.1$ in dimensions $d=5$ and $d = 10$
  seems to be an exception, where we see poor accuracy even for large $N$ and
  small $\delta$ (with the residuals for $N=10$ so large that they do not even
  appear in the plotting window). This may be due to the instability of autodiff
  in this case.}  
\label{fig:hj_complete}
\vspace*{-3in}
\end{figure}

\end{document}

%% file: ryantibs.tex
\usepackage{amssymb,amsmath,amsthm,bbm}
\usepackage{verbatim,float,url,dsfont}
\usepackage{graphicx,subcaption,psfrag}
\usepackage{algorithm,algorithmic}
\usepackage{mathtools,enumitem}
\usepackage{multirow}
\usepackage{ragged2e}
\usepackage{xr-hyper}
\usepackage{array}

\usepackage[colorlinks=true,citecolor=blue,urlcolor=blue,linkcolor=blue]{hyperref}
\usepackage[margin=1in]{geometry}
\usepackage[round]{natbib}

\usepackage[utf8]{inputenc} 
\usepackage[T1]{fontenc}    
\usepackage{booktabs}       
\usepackage{nicefrac}         
\usepackage{microtype}      

\ifdefined\TimesFont 
\usepackage{times} 
\fi

\ifdefined\ParSkip 
\usepackage{parskip} 
\fi

\newtheorem{theorem}{Theorem}

\newtheorem{corollary}{Corollary}

\theoremstyle{definition}

\newtheorem*{assumption*}{\assumptionnumber}
\providecommand{\assumptionnumber}{}


\makeatletter
\newcommand*\rel@kern[1]{\kern#1\dimexpr\macc@kerna}
\newcommand*\widebar[1]{%
  \begingroup
  \def\mathaccent##1##2{%
    \rel@kern{0.8}%
    \overline{\rel@kern{-0.8}\macc@nucleus\rel@kern{0.2}}%
    \rel@kern{-0.2}%
  }%
  \macc@depth\@ne
  \let\math@bgroup\@empty \let\math@egroup\macc@set@skewchar
  \mathsurround\z@ \frozen@everymath{\mathgroup\macc@group\relax}%
  \macc@set@skewchar\relax
  \let\mathaccentV\macc@nested@a
  \macc@nested@a\relax111{#1}%
  \endgroup
}
\makeatother

\DeclareMathOperator*{\argmin}{argmin}

\def\R{\mathbb{R}}

\def\E{\mathbb{E}}

\def\T{\mathsf{T}}

\def\cK{\mathcal{K}}
\def\cH{\mathcal{H}}

\def\cL{\mathcal{L}}

\def\cS{\mathcal{S}}

%% file: images/prox_point_csv/prox-point-plots-list.tex
\def\plotInfo{
    sphere/sphere_function,
    ellipsoidal/ellipsoidal_function,
    discus/discus_function,
    rosenbrock/rosenbrock_function,
    {sharp ridge}/sharp_ridge_function,
    weierstrass/weierstrass_function
}